\documentclass[11pt,reqno]{amsart}
\usepackage{amsmath,amssymb,amscd,amsthm}
\usepackage[pagewise]{lineno}

\usepackage{xcolor}
\usepackage{verbatim}
\usepackage{tikz-cd}
\usetikzlibrary{babel}
\usetikzlibrary{decorations.pathmorphing}
\usepackage[margin=1.3in]{geometry}
\usepackage{hyperref}

\hypersetup{hidelinks, colorlinks=true, allcolors=blue, pdfstartview=Fit, breaklinks=true, pdftitle={Approximate K-Conjugacy of Minimal Cantor-Circle Extensions I: The Rigid Case}, pdfauthor={Sihan Wei}}
\usepackage[all]{xy}
\DeclareMathOperator{\Ad}{Ad}
\usepackage{enumerate}
\usepackage{paralist}
\usepackage{amscd}
\IfFileExists{extarrows.sty}{\usepackage{extarrows}}{}

\newtheorem{thm}{Theorem}[subsection]
\newtheorem{lem}[thm]{Lemma}
\newtheorem{prop}[thm]{Proposition}

\newtheorem{cor}[thm]{Corollary}

\theoremstyle{definition}\newtheorem{df}[thm]{Definition}
\theoremstyle{definition}\newtheorem{rem}[thm]{Remark}

\theoremstyle{definition}

\renewcommand{\phi}{\varphi}

\newcommand{\N}{\mathbb{N}}
\newcommand{\Homeo}{\mathcal{H}}

\newcommand{\Z}{\mathbb{Z}}

\newcommand{\R}{\mathbb{R}}

\newcommand{\T}{\mathbb{T}}
\newcommand{\taf}{\tilde{\alpha}}
\newcommand{\tbt}{\tilde{\beta}}
\numberwithin{equation}{section}

\newcommand{\Aff}{\operatorname{Aff}}

\newcommand{\id}{\operatorname{id}}

\newcommand{\lr}{\longrightarrow}

\newcommand{\tr}{{\rm TR}}
\newcommand{\af}{{\alpha}}
\newcommand{\bt}{{\beta}}

\newcommand{\beq}{\begin{eqnarray}}
\newcommand{\eneq}{\end{eqnarray}}

\makeatletter

\def\@tocline#1#2#3#4#5#6#7{%
  \relax
  \ifnum #1>\c@tocdepth
  \else
    \par
    \addpenalty\@secpenalty
    \addvspace{#2}%
    \begingroup
      \hyphenpenalty\@M
      \@ifempty{#4}{%
        \@tempdima\csname r@tocindent\number#1\endcsname\relax
      }{%
        \@tempdima#4\relax
      }%
      \parindent\z@
      \leftskip#3\relax
      \rightskip\@pnumwidth
      \parfillskip-\@pnumwidth
      #5%
      \leavevmode
      #6%
      \nobreak
      \leaders\hbox{\normalfont.\kern0.5em}\hfill
      \nobreak
      \hbox to\@pnumwidth{\hss\normalfont #7}%
      \par
    \endgroup
  \fi
}

\renewcommand{\l@section}{%
  \@tocline{1}{0.8em}{0pt}{2.5em}{\bfseries}%
}

\renewcommand{\l@subsection}{%
  \@tocline{2}{0pt}{2.5em}{3.5em}{}%
}

\makeatother

\begin{document}

\title[Approximate $K$-Conjugacy of Cantor--Circle Extensions I]{Approximate $K$-Conjugacy of Minimal Cantor--Circle Extensions I: The Rigid Case}

\author{Huaxin Lin, Han Liu and Sihan Wei}
% \address{Department of Mathematics, ...}
% \email{...}
% \thanks{...}


\begin{abstract}
We study approximate $K$-conjugacy for minimal
Cantor--circle extensions with arbitrary circle
homeomorphism cocycles, allowing orientation
reversal. Under the assumption of classical
rigidity, we show that two such systems are
approximately $K$-conjugate if and only if
there exists an ordered unital graded
$K$-theory isomorphism satisfying the finite
coefficient $K$-condition. Our approach
combines dynamical corrections along
Kakutani--Rokhlin towers with the classification
and uniqueness theory of $C^*$-algebras.
Moreover, the crossed-product isomorphisms
implementing approximate $K$-conjugacy can
be chosen as inner perturbations of a fixed
isomorphism and its inverse.
\end{abstract}

\maketitle

\tableofcontents

\section{Introduction}\label{sec:introduction}

For a minimal homeomorphism $\alpha$ on a compact metric space $X$,
the crossed product $C(X)\rtimes_\alpha\mathbb Z$ provides a natural
connection between topological dynamics and the classification of
$C^*$-algebras. This connection is particularly transparent for Cantor
minimal systems. 
The crossed-product $C^*$-algebras associated
with Cantor minimal systems were studied
by Putnam \cite{Putnam1989}, who established
important structural connections with
AF-algebras and ordered $K$-theory.
Later on, Herman, Putnam and Skau studied their ordered dimension
groups by means of Bratteli diagrams \cite{HermanPutnamSkau1992}, while
Giordano, Putnam and Skau related the orbit structure of Cantor minimal
systems to ordered-group invariants and full groups
\cite{GiordanoPutnamSkau1999,GiordanoPutnamSkau2004,Putnam2010Survey}.
These developments suggest a general question: to what extent can
information obtained from the crossed product be used to compare the
underlying dynamical systems?

One way to approach this question is through approximate conjugacy.
Instead of requiring a single homeomorphism to conjugate two systems,
one allows a sequence of changes of coordinates under which the
dynamics become arbitrarily close. Approximate conjugacy is thus weaker
than conjugacy, but it still retains information about the homeomorphisms
rather than only their orbit partitions. Matui studied this relation
for Cantor minimal systems using topological full groups
\cite{Matui2005FullGroups}. Lin and Matui subsequently developed its
interaction with crossed-product $C^*$-algebras and introduced
approximate $K$-conjugacy as a way of incorporating the relevant
$K$-theoretic data \cite{LinMatui2005AC}. The essential point is that
one asks not only for approximate intertwining of the dynamics, but
also for compatibility between the same changes of coordinates and the
maps on the coefficient algebras induced by crossed-product
isomorphisms. In this sense, approximate $K$-conjugacy lies between a
purely dynamical approximation problem and the classification problem
for the associated $C^*$-algebras. 
The study of approximate conjugacy was
subsequently extended to other classes of
minimal dynamical systems. In particular,
Lin investigated approximate conjugacy
and approximate $K$-conjugacy for
Furstenberg transformations, further
illustrating the role of crossed-product
$K$-theory in distinguishing dynamical
systems \cite{LinFurstenberg2008}.

A natural setting beyond Cantor minimal systems is the product of a
Cantor space and the circle. Lin and Matui investigated minimal
systems on such spaces, including their $K$-theory, rigidity and
classification properties \cite{LinMatui2005,LinMatui2006}. The circle
factor introduces information which is absent in the zero-dimensional
case. In particular, a circle homeomorphism has degree in
$\{-1,1\}$, and this degree contributes to the action on the
$K_1$-group of the coefficient algebra. Consequently, an approximate
intertwining of the Cantor bases alone need not preserve the
$K$-theoretic information of the full extension. The distinction is
especially relevant when the circle maps are allowed to reverse
orientation.

In their study of minimal Cantor--circle systems, Lin and Matui
made extensive use of the structure of isometric extensions
\cite{LinMatui2005,LinMatui2006}. In that setting, the fiber dynamics
can be described by rotations and orientation data, so the geometry of
the cocycle is closely tied to its algebraic invariants. For general
circle-homeomorphism cocycles, this connection is no longer immediate.
Two fiber maps can have the same degree without being close as
homeomorphisms, and the degree records only the orientation of a map,
not the behavior of its iterates. Thus a comparison of the base systems
and their signed $K$-theoretic data does not by itself produce an
approximate conjugacy of the full extensions.

This leads to the question considered here: can approximate
$K$-conjugacy of minimal Cantor--circle extensions be characterized
without imposing an isometric assumption on the fiber maps? We work
with \emph{classically rigid} extensions, for which projection onto the
Cantor base induces a bijection between the spaces of invariant
probability measures. Our main result, Theorem~\ref{thm:main}, gives
an affirmative answer in this setting. It characterizes approximate
$K$-conjugacy by an ordered unital graded $K$-isomorphism of the
crossed products together with the \emph{finite coefficient
$K$-condition}. This condition requires that, on each prescribed
finitely generated subgroup of the $K$-groups of coefficient algebras, the maps
induced by coefficient isomorphisms eventually agree with the maps
obtained by composing the canonical embeddings with the given
$K$-isomorphism and its inverse. Moreover, the crossed-product
isomorphisms in the resulting approximation can be chosen as inner
perturbations of one fixed isomorphism and its inverse; in particular,
their $KL$-classes do not vary with the approximation.

The main difficulty is to pass from this finite algebraic condition
to an approximation by actual homeomorphisms. The $K$-groups of coefficient algebras carry both information from the Cantor base and the signs
of the circle maps. The latter enter through signed coinvariants, so
normalizing the base maps alone is insufficient: the orientation data
must also be made compatible. Even after this has been achieved,
equality of degrees does not control the fiber homeomorphisms in the
uniform topology. In particular, an adjustment that works at the
level of individual fiber maps need not respect the compositions
around a return tower.

To overcome this difficulty, we first translate the finite
coefficient $K$-condition into the required base and orientation
compatibilities. We then use clopen return towers to correct the
fiber dynamics. The correction is made along entire return
segments, with endpoint interpolation used to match the return
homeomorphisms, rather than merely their degrees. This construction
produces the desired dynamical approximation while retaining the
prescribed $K$-maps of coefficient algebras. The remaining issue is to realize
these maps through crossed-product isomorphisms.
Here classical rigidity supplies the necessary control of invariant
measures and traces. Together with the tracial-rank-zero results for
the crossed products and the relevant classification and uniqueness
theorems \cite{LinClassificationTR0,LinPhillips2010,MatuiHomomorphisms},
this allows the dynamical approximations to be implemented by inner
perturbations of a fixed isomorphism. The Pimsner--Voiculescu exact
sequence \cite{PimsnerVoiculescu1980} underlies the $K$-theoretic
computations used throughout.

This paper is the first part of our study of
approximate $K$-conjugacy for minimal
Cantor--circle extensions. In the sequel,
we will investigate the non-rigid case,
where the methods based on classical rigidity
are no longer directly available.

The paper is organized as follows. Section~\ref{sec:main-result}
introduces the systems and states the classification theorem, together
with a sketch of its proof. Section~\ref{sec:preliminaries} recalls
some basic facts about Cantor minimal systems and $K$-theory.
Section~\ref{sec:orientation} studies orientation signs and the
$K$-theory of the crossed products. In Section~\ref{sec:normalization}
we normalize the base and orientation data, and
Section~\ref{sec:towers} is devoted to the return-tower corrections.
Section~\ref{sec:rigid} establishes the trace and approximate unitary
equivalence results needed for the norm approximation. Finally,
Section~\ref{sec:proof} completes the proof of the main theorem.

\section*{Acknowledgements}

The second and third listed authors would like to express his sincere gratitude
to Huaxin for many stimulating discussions, his
continued encouragement, and, above all, his
shared perseverance and sustained efforts over
the years in helping the author understand
and address the problems considered in this paper.

\section*{Declaration of AI Use}

During the preparation of this manuscript, the authors
used ChatGPT (OpenAI) to assist with language editing,
clarity of exposition, and organization of the text.
All mathematical statements, arguments, and proofs
were reviewed and verified by the authors, who take
full responsibility for the content of this paper.

\section{The Main Result}
\label{sec:main-result}

\subsection{Minimal Cantor--circle extensions and rigidity}
For a compact metrizable space $(Z,d)$, let
$\Homeo(Z)$ denote the group of all homeomorphisms
of $Z$, equipped with the metric
\[
d_{\Homeo}(f,g)
=
\sup_{z\in Z}d(f(z),g(z))
+
\sup_{z\in Z}d(f^{-1}(z),g^{-1}(z)).
\]
This metric induces the compact-open topology
on $\Homeo(Z)$.
\medskip

A Cantor space is a nonempty compact metrizable
space without isolated points that has a basis
of clopen sets. 

\medskip

Let $\Omega_1$ and $\Omega_2$ be Cantor spaces, and put
\[
X=\Omega_1\times\T,
\qquad
Y=\Omega_2\times\T,
\]
where $\T=\R/\Z$. Although all Cantor spaces are homeomorphic, we use
different notations to distinguish the base spaces
of the two dynamical systems. Besides, we will use $\bar x$ to denote
the class of $x\in\R$ in $\T$ throughout the article.

\medskip

Let
\[
\af:\Omega_1\lr\Omega_1,
\qquad
\bt:\Omega_2\lr\Omega_2
\]
be homeomorphisms. We consider two homeomorphisms
\[
\taf:X\lr X,
\qquad
\tbt:Y\lr Y
\]
of the form
\[
\taf(\omega,t)
=
(\af(\omega),\varphi_\omega(t)),
\qquad
\tbt(\eta,t)
=
(\bt(\eta),\psi_\eta(t)).
\]
Here
\[
\omega\mapsto\varphi_\omega,
\qquad
\eta\mapsto\psi_\eta
\]
are continuous maps from $\Omega_1$ and $\Omega_2$,
respectively, into $\Homeo(\T)$. Every homeomorphism on the product of a Cantor space and the circle has this form, since the connected components of
$\Omega_i\times\T$ are precisely the circle fibers. We also write
$\taf=\af\times\varphi$
and
$\tbt=\bt\times\psi$
for these skew products.

\medskip

Throughout this paper, we assume that $\taf$ and $\tbt$
are minimal. Consequently, $\af$ and $\bt$ are minimal
homeomorphisms of the respective Cantor spaces. No restriction is imposed on the orientation of the
fiber homeomorphisms. In particular, the maps
$\varphi_\omega$ and $\psi_\eta$ need not be
rotations or isometries.

\medskip

Let
\[
A=C(X)\rtimes_{\taf}\Z,
\qquad
B=C(Y)\rtimes_{\tbt}\Z
\]
be the associated crossed product $C^*$-algebras.
We denote the canonical inclusions by
\[
j_{\taf}:C(X)\lr A,
\qquad
j_{\tbt}:C(Y)\lr B.
\]

For a homeomorphism $h$ of a compact metrizable space,
let $M(h)$ denote the compact convex space of
$h$-invariant Borel probability measures, equipped
with the weak$^*$ topology.

\begin{df}\label{def:classical-rigidity}
An extension $\taf=\af\times\phi$ on $X=\Omega\times\T$ is said to be
\emph{classically rigid} if the canonical projection
\[
\pi_X:X\lr\Omega
\]
induces a bijection
\[
(\pi_X)_*:M(\taf)\lr M(\af).
\]
\end{df}

\medskip
Since the spaces of invariant probability measures
are compact, such a bijection is automatically an affine
homeomorphism.
\medskip

\begin{rem}
    The notion of classical rigidity was introduced by Lin and Matui
in their earliest study of minimal dynamical systems on the
product of a Cantor space and the circle, see \cite[Definition~3.1]{LinMatui2005}.
They proved that a minimal Cantor--circle extension
is rigid if and only if its associated crossed product
has real rank zero.
Moreover, when the fiber maps are rotations,
rigidity implies that the crossed product has
tracial rank zero and is a unital simple
AT-algebra of real rank zero. Please refer to \cite[Corollary~3.11, Theorem~4.7, Theorem~5.6]{LinMatui2005} for more details.
\end{rem}

For a unital $C^*$-algebra $D$, let $T(D)$ denote
its tracial state space, and let
\[
\rho_D:K_0(D)\lr\Aff(T(D))
\]
denote the canonical pairing between $K_0(D)$
and traces.

\subsection{Approximate \texorpdfstring{$K$}{K}-conjugacy}
\begin{df}
Let $\taf\in\Homeo(X)$ and
$\tbt\in\Homeo(Y)$. We say that
$(X,\taf)$ and $(Y,\tbt)$ are
\emph{approximately conjugate} if there exist
sequences of homeomorphisms
\[
\sigma_n:X\lr Y,
\qquad
\gamma_n:Y\lr X
\]
such that
\[
\begin{aligned}
d_{\Homeo(Y)}
\bigl(\sigma_n\taf\sigma_n^{-1},\tbt\bigr)
&\lr 0,\\
d_{\Homeo(X)}
\bigl(\gamma_n\tbt\gamma_n^{-1},\taf\bigr)
&\lr 0.
\end{aligned}
\]
\end{df}

Lin and Matui
\cite[Proposition~3.2]{LinMatui2005AC}
showed that one-sided approximate conjugacy
induces a unital asymptotic morphism between
the associated crossed products. Let $u_{\taf}\in A$ and $u_{\tbt}\in B$
denote the canonical implementing unitaries, with covariance
\[
u_{\taf}j_{\taf}(f)u_{\taf}^*=j_{\taf}(f\circ\taf^{-1}),\qquad
u_{\tbt}j_{\tbt}(g)u_{\tbt}^*=j_{\tbt}(g\circ\tbt^{-1}).
\]
Let $E_A$ and $E_B$ be the canonical
conditional expectations onto $j_{\taf}(C(X))$
and $j_{\tbt}(C(Y))$, respectively.

\begin{prop}[{\cite[Proposition~3.2]{LinMatui2005AC}}]\label{prop:approx-asymptotic}
Suppose that $\taf$ and $\tbt$ are minimal.
Let $\sigma_n:X\to Y$ be homeomorphisms
such that
\[
d_{\Homeo(Y)}
\bigl(\sigma_n\taf\sigma_n^{-1},\tbt\bigr)
\longrightarrow 0.
\]
Then there exists a sequence of unital
completely positive maps
\[
\Psi_n:B\longrightarrow A
\]
which is asymptotically multiplicative and
satisfies
\[
\lim_{n\to\infty}
\|\Psi_n(u_{\tbt})-u_{\taf}\|=0
\]
and
\[
\lim_{n\to\infty}
\|\Psi_n(j_{\tbt}(g))
-j_{\taf}(g\circ\sigma_n)\|=0,
\qquad g\in C(Y).
\]
Moreover,
\[
\lim_{n\to\infty}
\|E_A(\Psi_n(b))-\Psi_n(E_B(b))\|=0
\]
for every $b\in B$.
\end{prop}

Proposition \ref{prop:approx-asymptotic} shows that approximate conjugacy
gives rise to asymptotic morphisms between the
associated crossed products. However, these
asymptotic morphisms need not be induced by
$*$-isomorphisms, nor do they necessarily
determine fixed classes in $KL$-theory.
This motivates a stronger notion of approximate
conjugacy, in which the dynamical approximations
are required to be compatible with sequences
of $*$-isomorphisms having fixed inverse
$KL$-classes.

\medskip

\begin{df}\label{def:approx-K-conjugacy}
The systems $(X,\taf)$ and $(Y,\tbt)$ are said to be
\emph{approximately $K$-conjugate} if there exist
sequences of homeomorphisms
\[
\sigma_n:X\lr Y,
\qquad
\gamma_n:Y\lr X,
\]
implementing the approximate conjugacy between $(X,\taf)$ and $(Y,\tbt)$, and also sequences of $*$-isomorphisms
\[
\Phi_n:A\lr B,
\qquad
\Psi_n:B\lr A
\]
with an invertible element
\[
\kappa\in KL(A,B),
\]
such that
\[
[\Phi_n]=\kappa,
\qquad
[\Psi_n]=\kappa^{-1}
\]
for every $n\in\N$, and the following
conditions hold:
\begin{align}
\big\|
\Phi_n(j_{\taf}(f))
-j_{\tbt}(f\circ\gamma_n)
\big\|
&\lr 0,
\qquad \forall f\in C(X),
\label{eq:approx-coeff-forward}
\\
\big\|
\Psi_n(j_{\tbt}(g))
-j_{\taf}(g\circ\sigma_n)
\big\|
&\lr 0,
\qquad \forall g\in C(Y).
\label{eq:approx-coeff-backward}
\end{align}
\end{df}

\medskip

\begin{rem}
The same homeomorphisms $\sigma_n$ and $\gamma_n$
occur in both the dynamical approximation and
the coefficient-algebra approximation. Besides, we do not require
\[
\gamma_n=\sigma_n^{-1}
\qquad\text{or}\qquad
\Psi_n=\Phi_n^{-1}.
\]
Moreover, the trace maps induced by $\Phi_n$
and $\Psi_n$ need not be independent of $n$.

Approximate $K$-conjugacy is stronger than
approximate conjugacy. In addition to the
dynamical approximation, it further requires the
corresponding homeomorphisms to be compatible,
asymptotically, with isomorphisms of the
associated crossed products having fixed
inverse $KL$-classes.
\end{rem}

\subsection{The classification theorem}

We now give a characterization of approximate
$K$-conjugacy for classically rigid minimal
Cantor--circle extensions in terms of the
ordered $K$-theory of their crossed products
and the finite coefficient $K$-condition.

\begin{thm}\label{thm:main}
Let $\taf:X\lr X$ and $\tbt:Y\lr Y$
be minimal Cantor--circle extensions as above.
Suppose that both $\taf$ and $\tbt$ are
classically rigid.

Then the following conditions are equivalent.

\begin{enumerate}
\item
The systems $(X,\taf)$ and $(Y,\tbt)$
are approximately $K$-conjugate.

\item
There exist an ordered unital graded isomorphism of $K$-theory
\[
k=(k_0,k_1):
\big(
K_0(A),K_0(A)_+,[1_A],K_1(A)
\big)
\lr
\big(
K_0(B),K_0(B)_+,[1_B],K_1(B)
\big)
\]
and sequences of unital $*$-isomorphisms
\[
\chi_n:C(X)\lr C(Y),
\qquad
\eta_n:C(Y)\lr C(X),
\]
satisfying the following finite coefficient
$K$-condition:
\medskip

For every choice of finitely generated subgroups
\[
H_i\subset K_i(C(X)),
\qquad
J_i\subset K_i(C(Y)),
\qquad i=0,1,
\]
there exists an integer $N\geq 1$ such that,
for every $n\geq N$ and both $i=0,1$,
\begin{align}
\left.
k_i\circ(j_{\taf})_{*i}
\right|_{H_i}
&=
\left.
(j_{\tbt}\circ\chi_n)_{*i}
\right|_{H_i},
\label{eq:main-K-forward}
\\
\left.
k_i^{-1}\circ(j_{\tbt})_{*i}
\right|_{J_i}
&=
\left.
(j_{\taf}\circ\eta_n)_{*i}
\right|_{J_i}.
\label{eq:main-K-backward}
\end{align}
\end{enumerate}
\medskip

Moreover, whenever these conditions hold,
the crossed-product isomorphisms witnessing
approximate $K$-conjugacy may be chosen in
the form
\[
\Phi_n=\Ad(u_n)\circ\Phi,
\qquad
\Psi_n=\Ad(w_n)\circ\Phi^{-1},
\]
where $\Phi:A\lr B$ is a fixed
$*$-isomorphism and
\[
u_n\in U(B),
\qquad
w_n\in U(A)
\]
are unitaries. In particular,
\[
[\Phi_n]=[\Phi]\in KL(A,B),
\qquad
[\Psi_n]=[\Phi]^{-1}\in KL(B,A)
\]
for every $n\in\N$.
\end{thm}

\medskip

\begin{rem}
In the finite coefficient $K$-condition,
the integer $N$ may depend on the four
prescribed finitely generated subgroups,
but is common to both directions and
both $K$-theory degrees.
\end{rem}

\subsection{Sketch of the proof}
\label{subsec:proof-sketch}

We briefly outline the proof of
Theorem~\ref{thm:main}. The necessity follows
from the coefficient norm approximations and
the fixed inverse $KL$-classes. The main
difficulty is to establish the converse.

We begin by using the ordinary and signed
coinvariant groups, together with corrections
from the topological full groups, to normalize
the coefficient homeomorphisms so that their
base maps approximately intertwine the Cantor
minimal systems and their fiber degrees agree
with the corresponding target degrees, without
changing the prescribed coefficient $K$-maps.

We then use Kakutani--Rokhlin partitions and
return-tower interpolation to correct the fiber
dynamics, since agreement of fiber degrees
alone does not imply uniform closeness of
circle homeomorphisms. These corrections are
isotopic to the identity and therefore preserve
the coefficient $K$-theory data. This produces
approximate conjugacies in both directions.

Finally, classical rigidity provides tracial
rank zero for the crossed products and uniform
trace compatibility for the coefficient
embeddings. Combining these properties with
the finite coefficient $K$-condition, we apply
the uniqueness theorem for embeddings to
obtain the required coefficient norm
approximations. The implementing
$*$-isomorphisms are inner perturbations of
a fixed isomorphism and its inverse, so their
$KL$-classes remain unchanged.

\section{Preliminaries}\label{sec:preliminaries}

In this section, we introduce some basic facts
about dynamical systems, $K$-theory,
and crossed products that will be used in
the proof in later sections.
\medskip

\subsection{Basic notions from topological dynamics}

We begin by recalling some standard notions
from topological dynamics.

Let $Z$ be a compact metrizable space and
let $\af:Z\to Z$ be a homeomorphism.
The pair $(Z,\af)$ is called a
\emph{topological dynamical system}.

For $z\in Z$, the orbit of $z$ under $\af$ is
defined by
\[
\operatorname{Orb}_\af(z)
=
\{\af^n(z):n\in\Z\}.
\]

A subset $F\subset Z$ is said to be
$\af$-invariant if $\af(F)=F$.

The system $(Z,\af)$ is called \emph{minimal}
if $Z$ contains no nonempty proper closed
$\af$-invariant subset. Equivalently,
$(Z,\af)$ is minimal if every orbit is dense
in $Z$, that is,
\[
\overline{\operatorname{Orb}_\af(z)}=Z
\]
for every $z\in Z$.

A minimal topological dynamical system $(\Omega,\af)$ on a Cantor space $\Omega$ is
called a \emph{Cantor minimal system}.

\subsection{The Kakutani-Rokhlin partition}
We now recall the construction of
Kakutani--Rokhlin partitions for Cantor
minimal systems; see
\cite{HermanPutnamSkau1992}.

Let $(\Omega,\alpha)$ be a Cantor minimal
system, and let $U\subset\Omega$ be a
nonempty clopen subset. For each $\omega\in U$, define the first
return time to $U$ by
\[
r_U(\omega)
=
\min\{n\geq 1:\alpha^n(\omega)\in U\}.
\]
Since $\alpha$ is minimal and $U$ is
clopen, the function
\[
r_U:U\longrightarrow\N
\]
is well-defined, continuous, and takes
only finitely many values.

Let $h_1,\ldots,h_m$ be the distinct values
of $r_U$, and put
\[
U_k=\{\omega\in U:r_U(\omega)=h_k\},
\qquad 1\leq k\leq m.
\]
Then each $U_k$ is clopen, and
\[
U=\bigsqcup_{k=1}^m U_k.
\]

For each $k$, the collection
\[
\mathcal T_k
=
\{U_k,\alpha(U_k),\ldots,
\alpha^{h_k-1}(U_k)\}
\]
is called a \emph{Kakutani--Rokhlin tower}
of height $h_k$, with base $U_k$ and
roof $\alpha^{h_k-1}(U_k)$.

The collection
\[
\mathcal P_U
=
\big\{
\alpha^j(U_k):
1\leq k\leq m,\ 0\leq j<h_k
\big\}
\]
forms a finite clopen partition of
$\Omega$, i.e.,
\[
\Omega
=
\bigsqcup_{k=1}^m
\bigsqcup_{j=0}^{h_k-1}
\alpha^j(U_k),
\]
called the
\emph{Kakutani--Rokhlin partition}
associated with $U$. We refer to the sets $\alpha^j(U_k)$
as the levels of the corresponding tower.

\begin{rem}
For every $N\geq 1$, there exists
a nonempty clopen set $U\subset\Omega$
such that every tower in the Kakutani--Rokhlin partition associated
with $U$ has height at least $N$. 

In fact, this could be done by choosing sufficiently small $U$ so that $U,\af(U),\cdots,\af^{N-1}(U)$ are pairwise disjoint since $(\Omega,\af)$ has no periodic points.
\end{rem}

\subsection{Topological full group}
We recall the notion of the topological
full group of a Cantor minimal system;
see \cite{GiordanoPutnamSkau1999,
Matui2006FullGroups}.

The \emph{topological full group} of a
Cantor minimal system $(\Omega,\alpha)$, denoted by $[[\alpha]]$, is the
group of all homeomorphisms
$\sigma:\Omega\to\Omega$ for which there
exists a continuous function
$n_\sigma:\Omega\to\Z$ such that
\[
\sigma(\omega)
=
\alpha^{n_\sigma(\omega)}(\omega),
\qquad \omega\in\Omega.
\]
Since $\Omega$ is compact and totally
disconnected, the function $n_\sigma$
takes only finitely many values.
Consequently, there exists a finite
clopen partition
\[
\Omega=\bigsqcup_{k=1}^m U_k
\]
and integers $n_1,\ldots,n_m$ such that
\[
\sigma|_{U_k}
=
\alpha^{n_k}|_{U_k},
\qquad 1\leq k\leq m.
\]
In particular, every element of
$[[\alpha]]$ preserves each
$\alpha$-orbit setwise.

\subsection{\texorpdfstring{$K$}{K}-theory of Cantor-circle spaces}
We recall the standard $K$-theory
computations for Cantor--circle spaces;
see, for example,
\cite{Blackadar1998,WeggeOlsen1993}
for the relevant facts from operator
$K$-theory.

Let $\Omega$ be a Cantor space. Since $\Omega$
is totally disconnected, we have
\[
K_0(C(\Omega))\cong C(\Omega,\Z),
\qquad
K_1(C(\Omega))=0.
\]
By the K\"unneth theorem, together with
\[
K_0(C(\T))\cong\Z,
\qquad
K_1(C(\T))\cong\Z,
\]
we obtain
\[
K_i(C(\Omega\times\T))
\cong C(\Omega,\Z),
\qquad i=0,1.
\]

Under these identifications, the positive cone
of $K_0(C(\Omega\times\T))$ is $\{f\ge0: f\in C(\Omega,\Z)\}$, with the order unit being the constant function $1$. 
These identifications admit explicit
descriptions in terms of fiberwise rank
and winding number in the following manner, which will be useful
in the sequel.

\medskip

For a projection
$p\in M_m(C(\Omega\times\T))$,
the rank of $p(\omega,t)$ is independent
of $t$. Thus the function
\[
\omega\longmapsto
\operatorname{rank}(p(\omega,t))
\]
is a continuous integer-valued function
on $\Omega$, called the fiberwise rank function of $p$.

Similarly, for a unitary
$u\in U_m(C(\Omega\times\T))$, the function
\[
\omega\longmapsto
\operatorname{wind}
\big(\det u(\omega,\cdot)\big)
\]
is continuous and integer-valued, where
$\operatorname{wind}$ denotes the winding
number of a circle-valued map. We refer to it as the fiberwise winding number of $u$.

\subsection{The Pimsner--Voiculescu exact sequence}\label{subsec:PV}

We recall the Pimsner--Voiculescu exact
sequence for crossed products of $C^*$-algebras by $\Z$ (see \cite{PimsnerVoiculescu1980}).
\medskip

Let $D$ be a unital $C^*$-algebra and
$\theta\in\operatorname{Aut}(D)$.
Write
\[
D\rtimes_\theta\Z
\]
for the crossed product, with implementing unitary $v$ satisfying
\[
vdv^*=\theta(d),
\qquad d\in D.
\]

The Pimsner--Voiculescu six-term exact
sequence is
\[
\begin{CD}
K_0(D) @>{\mathrm{id}-\theta_{*0}}>>
K_0(D) @>{j_{*0}}>>
K_0(D\rtimes_\theta\Z)\\
@A{\partial_1}AA
&& @VV{\partial_0}V\\
K_1(D\rtimes_\theta\Z)
@<{j_{*1}}<< K_1(D)
@<{\mathrm{id}-\theta_{*1}}<< K_1(D).
\end{CD}
\]
where $j:D\hookrightarrow D\rtimes_\theta\Z$ is the canonical
inclusion. We adopt the sign convention
\[
\partial_1([v])=[1_D].
\]
In particular, exactness gives
\[
\ker j_{*i}
=
\operatorname{im}(\mathrm{id}-\theta_{*i}),
\qquad i=0,1.
\]
We shall use these identities to identify
the images of the canonical coefficient
maps in the $K$-groups of the crossed products.
The corresponding calculations for
Cantor--circle extensions will be carried
out in the next section.

\section{Orientation and \texorpdfstring{$K$}{K}-theory}\label{sec:orientation}

\subsection{Orientation signs}

We first recall the degree of a continuous
map of the circle. Let $h:\T\to\T$ be
continuous, and let $\widehat h:\R\to\R$
be a lift of $h$ under the quotient map
$\R\to\T=\R/\Z$. The degree of $h$ is
defined by
\[
\deg(h)=\widehat h(r+1)-\widehat h(r),
\qquad r\in\R.
\]
This is an integer independent of $r$
and of the choice of lift.

If $h$ is a homeomorphism, then
$\deg(h)\in\{-1,1\}$. We say that $h$
is {\it orientation-preserving} if
$\deg(h)=1$, and {\it orientation-reversing}
if $\deg(h)=-1$.

For the Cantor--circle extensions
$\widetilde\alpha=\alpha\times\varphi$
and $\widetilde\beta=\beta\times\psi$
introduced in Section~\ref{sec:main-result}, define their
orientation signs by
\[
s_{\widetilde\alpha}(\omega)
=\deg(\varphi_\omega),
\qquad \omega\in\Omega_1,
\]
and
\[
s_{\widetilde\beta}(\eta)
=\deg(\psi_\eta),
\qquad \eta\in\Omega_2.
\]
Since degree is locally constant on
$\Homeo(\T)$, both
\[
s_{\widetilde\alpha}:
\Omega_1\to\{-1,1\},
\qquad
s_{\widetilde\beta}:
\Omega_2\to\{-1,1\}
\]
are continuous.

\medskip

\begin{df}\label{def:sign-coboundary}
Let $(\Omega,\alpha)$ be a Cantor minimal
system, and let
$s:\Omega\to\{-1,1\}$ be continuous.

We say that $s$ is a
\emph{sign coboundary} over $\alpha$
if there exists a continuous function
$a:\Omega\to\{-1,1\}$ such that
\[
s(\omega)=a(\alpha(\omega))a(\omega),
\qquad \omega\in\Omega.
\]

Two continuous sign functions
$s,s':\Omega\to\{-1,1\}$ are said to
have the same
\emph{continuous orientation class}
if there exists a continuous function
$a:\Omega\to\{-1,1\}$ such that
\[
s'(\omega)
=
a(\alpha(\omega))s(\omega)a(\omega),
\qquad \omega\in\Omega.
\]
The continuous orientation class of
$s$ is called \emph{trivial} if $s$
is a sign coboundary.

The {\it sign extension} associated to $\taf$ is a homeomorphism
\[\widehat{\taf}: \Omega\times\{-1,1\}\times\T\to \Omega\times\{-1,1\}\times\T,\]
sending $(\omega,e,z)$ to $(\af(\omega),s(\omega)e,\phi_\omega(z))$, where $s=s_{\taf}$ is the orientation sign.
\end{df}

\begin{lem}\label{lem:reflection-coordinates}
The sign $s_{\taf}$ is a sign coboundary if and only if there is a homeomorphism $C: X\to X$ fixing the first coordinate such that 
\[\taf^+=C\taf C^{-1}\]
has orientation-preserving fiber maps everywhere.

Besides, the conjugacy preserves minimality and
classical rigidity. The coefficient isomorphism
$\theta_a(f)=f\circ C^{-1}$ extends to a crossed-product isomorphism
$\Theta_a:A\to A^+$ with
$\Theta_a j_{\taf}=j_{\taf^+}\theta_a$.
\end{lem}

\begin{proof}
Suppose that $C$ has the form
\[C(\omega,z)=(\omega, c_\omega(z))\]
where $c_\omega\in\Homeo(\T)$. Define $a(\omega)={\rm deg}(c_\omega)$. Then $a:\Omega\to\{-1,1\}$ is continuous. The fiber map of $C\taf C^{-1}$ at $\omega$ is $c_{\af(\omega)}\varphi_\omega c_\omega^{-1}$, whose degree is
\[a(\af(\omega))s_{\taf}(\omega)a(\omega)=1\]
as assumed. This gives us 
\[s_{\taf}(\omega)=a(\af(\omega))a(\omega),\]
since $a^2=1$. Therefore $s_{\taf}$ is a sign coboundary.

Conversely, assume that $s_{\taf}$ is a sign coboundary. Write 
\[s_{\taf}(\omega)=a(\af(\omega))a(\omega)\]
for some continuous function $a:\Omega\to\{-1,1\}$. Define 
\[C_a(\omega,z)=(\omega,a(\omega)z).\]
A straightforward calculation shows that
\[C_a\taf C_a^{-1}(\omega,z)=(\af(\omega), a(\af(\omega))\phi_\omega(a(\omega)z)).\]
The degree of the resulting homeomorphism equals
\[a(\af(\omega))\cdot s_{\taf}(\omega)\cdot a(\omega)=[a(\af(\omega))]^2\cdot [a(\omega)]^2=1.\]
It is immediate that conjugacy preserves minimality, and since $C_a$ fixes the base projection,
it identifies the invariant lift of each base probability in the
two systems. 

Finally, the covariance extends $\theta_a$, sending the implementing
unitary to that for $\taf^+$. The inverse coordinate change gives its
inverse crossed-product map.
\end{proof}

\begin{lem}\label{lem:orientation-cover}
Let $\taf$ be minimal and $s=s_{\taf}$ not a sign coboundary.
Then both
\[
\widehat\alpha(\omega,e)=(\alpha(\omega),s(\omega)e),\qquad
\widehat{\taf}(\omega,e,z)
=(\alpha(\omega),s(\omega)e,\varphi_\omega(z))
\]
are minimal. The coordinate change $(\omega,e,z)\mapsto(\omega,e,ez)$
conjugates $\widehat{\taf}$ to an extension whose fiber maps all preserve orientation.
\end{lem}
\begin{proof}
    The proof follows from a straightforward computation and is therefore omitted.
\end{proof}

\begin{rem}\label{rem:cover-role}
The orientation cover will be used for local signed corrections and
for endpoint width. Its minimality is needed; no rigidity assertion
about this cover is used.
\end{rem}

\subsection{Signed coinvariants}\label{subsec:signed-coinvariants}

\begin{df}
Let $(\Omega,\alpha)$ be a Cantor minimal
system, and let
$s:\Omega\to\{-1,1\}$ be continuous.
Put
\[
M=C(\Omega,\Z).
\]
Define automorphisms $U_1$ and $U_s$
of the additive group $M$ by
\[
(U_1f)(\omega)
=f(\alpha^{-1}(\omega))
\]
and
\[
(U_sf)(\omega)
=s(\alpha^{-1}(\omega))
 f(\alpha^{-1}(\omega)),
\qquad f\in M.
\]
The groups
\[
D_\alpha=M/(I-U_1)M
\]
and
\[
Q_s=M/(I-U_s)M
\]
are called the \emph{ordinary coinvariant
group} and the \emph{signed coinvariant
group}, respectively.
\end{df}

\begin{lem}\label{lem1}
    Let $s:\Omega\to\{-1,1\}$ be continuous. If $s$ is a sign coboundary, then
    \[Q_s\cong D_\af.\]
\end{lem}
\begin{proof}
    Suppose that $s$ is a sign coboundary. Then we can choose a continuous function $a:\Omega\to\{-1,1\}$ satisfying
    \[s(\omega)=a(\af(\omega))a(\omega).\]
    Define an automorphism $B_a:M\to M$ by
    \[(B_af)(\omega)=a(\omega)f(\omega),\ \ f\in M.\]
    One checks that $B_a U_s=U_1 B_a$. Therefore $B_a$ induces an isomorphism from $Q_s$ to $D_\af$, sending $[f]\in Q_s$ to $[B_af]\in D_\af$.  This gives the isomorphism as desired.
\end{proof}

\medskip

\begin{prop}\label{prop:signed-torsion}
Let $(\Omega,\alpha)$ be a Cantor minimal
system, and let
$s:\Omega\to\{-1,1\}$ be continuous.
Let $Q_s$ be the signed coinvariant group
defined above.

(1) If $s$ is a sign coboundary, then $Q_s$
is torsion-free.

(2) If $s$ is not a sign coboundary, then
the torsion subgroup of $Q_s$ is cyclic
of order two, generated by $[v_s]$, where
\[
v_s(\omega)
=
\frac{1-s(\alpha^{-1}(\omega))}{2},
\qquad \omega\in\Omega.
\]
\end{prop}

\begin{proof}
We first consider the case {\bf where $s$ is
a sign coboundary}.

Suppose initially that $s\equiv 1$, in which case $Q_s$ is just the group $D_\af$. 
Let $n\geq 1$ and suppose that
$n[f]=0$ in $D_\alpha$. Then there exists
$g\in M$ such that
\[
nf=g-g\circ \af^{-1}.
\]
Reducing modulo $n$, we obtain
\[
g(\alpha(\omega))
\equiv g(\omega)\pmod n,
\qquad \omega\in\Omega.
\]
Thus the residue of $g$ modulo $n$ is
a continuous $\alpha$-invariant function.
By minimality of $\af$, it must be constant. 

Consequently, there exist an integer 
$c\in\Z$ and a function $g'\in C(\Omega,\Z)$ such that
\[
g=c+ng'.
\]
It follows that
\[
nf=n(g'-g'\circ\af^{-1}),
\]
and hence
\[
f=g'-g'\circ\af^{-1}=g'-U_1g'\in(I-U_1)M.
\]
Therefore $[f]=0$, showing that
$D_\alpha$ is torsion-free.

Now suppose that $s$ is a sign
coboundary. According to Lemma \ref{lem1}, we have
\[
Q_s\cong D_\alpha.
\]
Since $D_\alpha$ is torsion-free,
so is $Q_s$.

\medskip

We now assume that {\bf $s$ is not a
sign coboundary}.

\medskip

Let $[f]\in Q_s$ be a torsion element.
Then there exist $n\geq 1$ and
$g\in M$ such that
\[
nf=g-U_sg.
\]
Let
\[
r:\Omega\to\Z/n\Z
\]
be the residue of $g$ modulo $n$. It is immediate that $r$ is continuous. 
By evaluating the equality above 
at $\alpha(\omega)$ modulo $n$, we have
\[0\equiv nf(\af(\omega))\equiv g(\af(\omega))-s(\omega)g(\omega)\pmod n,\]
or in other words,
\[r(\af(\omega))=s(\omega)r(\omega)\ \ {\rm in}\ \Z/n\Z.\]

Consider the image of $r$ on any single orbit of $\af$. Since $s$ takes value in $\{-1,1\}$, this image consists of at most two elements, say, $-c$ and $c$. Moreover, since the homeomorphism $\af$ is minimal and $r$ is continuous, any $\af$-orbit is dense in $\Omega$, which means that the image of $r$ on the whole space $\Omega$ is again a subset of $\{-c,c\}$.

Now we already have a constant $c\in\Z/n\Z$ such that
\[
r(\Omega)\subseteq\{c,-c\}.
\]

Suppose that $c\neq -c$.
Define
\[
a(\omega)=
\begin{cases}
1,&r(\omega)=c,\\
-1,&r(\omega)=-c.
\end{cases}
\]
Then $a$ is continuous and satisfies
\[
a(\alpha(\omega))
=
s(\omega)a(\omega).
\]
This would imply that $s$ is a
sign coboundary, a contradiction.

Therefore $c=-c$, and $r$ is constant.
In particular, either $r=0$, or
$n$ is even and $r=n/2$ in
$\Z/n\Z$.

If $r=0$, we may write $g=ng'$
for some $g'\in M$. Hence
\[
f=g'-U_sg',
\]
so $[f]=0$.

In the remaining case, $n$ is even and we may write
\[
g=ng'+\frac{n}{2}
\]
for some $g'\in M$. Consequently,
\[
\begin{aligned}
f
&=(g'-U_sg')+
\frac{1-s\circ\af^{-1}}{2}\\
&=g'-U_sg'+v_s.
\end{aligned}
\]
Thus $[f]=[v_s]$. Moreover,
\[
2v_s=1-s\circ\af^{-1}=1_\Omega-U_s1_\Omega\in(I-U_s)M,
\]
which shows that
\[
2[v_s]=0
\quad\text{in }Q_s.
\]
\medskip

It remains to show that
$[v_s]\neq 0$. Suppose, to the contrary, that
\[
v_s=k-U_sk
\]
for some $k\in M$. Since $2v_s=1_\Omega-U_s1_\Omega$, 
we obtain
\[
U_s(1_\Omega-2k)=1_\Omega-2k.
\]
Put
\[
w=1_\Omega-2k.
\]
Then
\[
w(\alpha(\omega))
=
s(\omega)w(\omega).
\]
In particular, $|w|$ is a continuous
$\alpha$-invariant function.
By minimality, $|w|$ is constant,
say equal to $m$.

Since $w$ is odd-valued, $m$ is
a positive odd integer. Therefore
\[
a=\frac{w}{m}
\]
is a continuous function from
$\Omega$ to $\{-1,1\}$ satisfying
\[
a(\alpha(\omega))
=
s(\omega)a(\omega).
\]
This again implies that $s$ is a
sign coboundary, a contradiction. Hence $[v_s]\neq 0$ and we conclude that
\[
{\rm Tor}(Q_s)
=
\{0,[v_s]\}
\cong\Z/2\Z.
\]
\end{proof}

\subsection{K-theory of the crossed products}\label{subsec:crossed-product-K}

We now apply the Pimsner--Voiculescu exact
sequence to the Cantor--circle extensions. In the following Proposition, by a slight abuse of notation,
$(j_{\widetilde\alpha})_{*i}$ also denotes
the injective homomorphism induced on
the corresponding coinvariant quotient
by the canonical inclusion
$j_{\widetilde\alpha}:C(X)\to A$.

\begin{prop}\label{prop:coefficient-exact-sequences}
Let $\widetilde\alpha=\alpha\times\varphi$
be a minimal Cantor--circle extension on
$X=\Omega_1\times\T$, and put
\[
A=C(X)\rtimes_{\widetilde\alpha}\Z.
\]
Let $s=s_{\widetilde\alpha}$, and let
$D_\alpha$ and $Q_s$ be the coinvariant
groups defined in Subsection~\ref{subsec:signed-coinvariants}.

Then there are short exact sequences
\[
0\longrightarrow D_\alpha
\xrightarrow{\ (j_{\widetilde\alpha})_{*0}\ }
K_0(A)
\xrightarrow{\ \partial_0\ }
\ker(I-U_s)
\longrightarrow 0
\]
and
\[
0\longrightarrow Q_s
\xrightarrow{\ (j_{\widetilde\alpha})_{*1}\ }
K_1(A)
\xrightarrow{\ \partial_1\ }
\Z
\longrightarrow 0.
\]

Moreover, the second short exact sequence
splits, i.e.,
\[
K_1(A)\cong Q_s\oplus\Z,
\]
where the summand $Q_s$ is identified
with the image of
$(j_{\widetilde\alpha})_{*1}$.
\end{prop}

\begin{proof}
Write
\[
\theta:C(X)\longrightarrow C(X),
\qquad
\theta(f)=f\circ\widetilde\alpha^{-1}.
\]
We first identify the induced maps
on the $K$-groups of $C(X)$.

Under the identification
\[
K_0(C(X))\cong C(\Omega_1,\Z),
\]
the automorphism $\theta$ induces
\[
(\theta_{*0}f)(\omega)
=f(\alpha^{-1}(\omega))
=(U_1f)(\omega).
\]

Similarly, under the identification
\[
K_1(C(X))\cong C(\Omega_1,\Z),
\]
we have
\[
\begin{aligned}
(\theta_{*1}f)(\omega)
&=
s(\alpha^{-1}(\omega))
f(\alpha^{-1}(\omega))\\
&=(U_sf)(\omega).
\end{aligned}
\]
Indeed, the restriction of
$\widetilde\alpha^{-1}$ to the fiber
over $\omega$ is determined by
$\varphi_{\alpha^{-1}(\omega)}^{-1}$,
whose degree equals
$s(\alpha^{-1}(\omega))$.
Therefore
\[
\theta_{*0}=U_1,
\qquad
\theta_{*1}=U_s.
\]

Applying the Pimsner--Voiculescu exact
sequence,
we obtain
\[
\ker(j_{\widetilde\alpha})_{*0}
=(I-U_1)C(\Omega_1,\Z)
\]
and
\[
\ker(j_{\widetilde\alpha})_{*1}
=(I-U_s)C(\Omega_1,\Z).
\]
Thus the coefficient maps induce
injective homomorphisms
\[
(j_{\widetilde\alpha})_{*0}:
D_\alpha\longrightarrow K_0(A)
\]
and
\[
(j_{\widetilde\alpha})_{*1}:
Q_s\longrightarrow K_1(A).
\]
Exactness also gives
\[
K_0(A)/
\operatorname{im}(j_{\widetilde\alpha})_{*0}
\cong\ker(I-U_s)
\]
and
\[
K_1(A)/
\operatorname{im}(j_{\widetilde\alpha})_{*1}
\cong\ker(I-U_1).
\]

Since $\alpha$ is minimal, every
continuous integer-valued function
invariant under $\alpha$ is constant.
Hence
\[
\ker(I-U_1)\cong\Z.
\]
This establishes the two short exact
sequences.

Finally, let $u_{\widetilde\alpha}\in A$
be the canonical implementing unitary.
By the convention fixed in Subsection~\ref{subsec:PV},
\[
\partial_1([u_{\widetilde\alpha}])
=[1_{C(X)}].
\]
Under the identification
$\ker(I-U_1)\cong\Z$, the latter
corresponds to $1\in\Z$.
Consequently, the homomorphism
\[
\Z\longrightarrow K_1(A),
\qquad
n\longmapsto n[u_{\widetilde\alpha}]
\]
is a right inverse of $\partial_1$.

Therefore the second short exact
sequence splits, and
\[
K_1(A)\cong Q_s\oplus\Z,
\]
where $Q_s$ is identified with
$\operatorname{im}(j_{\widetilde\alpha})_{*1}$.
\end{proof}

\begin{prop}\label{prop:orientation-dichotomy}
Keep the notations as in Proposition \ref{prop:coefficient-exact-sequences}. The following statements hold.

\begin{enumerate}
\item If $s$ is a sign coboundary, then
\[
\ker(I-U_s)\cong\Z,
\qquad
K_0(A)\cong D_\alpha\oplus\Z,
\]
and $K_1(A)$ is torsion-free.

\item If $s$ is not a sign coboundary, then
\[
\ker(I-U_s)=0,
\qquad
K_0(A)\cong D_\alpha,
\]
and
\[
\operatorname{Tor}(K_1(A))
\cong\Z/2\Z.
\]
More precisely, the nonzero torsion
element of $K_1(A)$ is
\[
(j_{\widetilde\alpha})_{*1}([v_s]),
\]
where $[v_s]\in Q_s$ is the element
defined in Proposition \ref{prop:signed-torsion}.
\end{enumerate}
\end{prop}

\begin{proof}
We first determine $\ker(I-U_s)$. Suppose that
$f\in\ker(I-U_s)$. Then
\[
f(\omega)
=s(\alpha^{-1}(\omega))
f(\alpha^{-1}(\omega)),
\]
or equivalently,
\[
f(\alpha(\omega))
=s(\omega)f(\omega),
\qquad \omega\in\Omega_1.
\]
It follows that
\[
|f(\alpha(\omega))|=|f(\omega)|.
\]
Since $\alpha$ is minimal and $|f|$
is continuous, $|f|$ is constant
on $\Omega_1$.

If $f\ne0$, there exists an integer
$m\ge1$ such that $|f|=m$.
Define
\[
a(\omega)=\frac{f(\omega)}{m},
\qquad \omega\in\Omega_1.
\]
Then $a:\Omega_1\to\{-1,1\}$ is
continuous and satisfies
\[
a(\alpha(\omega))=s(\omega)a(\omega).
\]
Consequently,
\[
s(\omega)=a(\alpha(\omega))a(\omega),
\]
showing that $s$ is a sign coboundary. Thus, if $s$ is not a sign coboundary,
we must have
\[
\ker(I-U_s)=0.
\]

Now suppose that $s$ is a sign
coboundary. Choose a continuous
function $a:\Omega_1\to\{-1,1\}$
such that
\[
s(\omega)=a(\alpha(\omega))a(\omega).
\]
As observed in Lemma \ref{lem1},
multiplication by $a$ defines an
automorphism $B_a$ of
$C(\Omega_1,\Z)$ satisfying
\[
B_aU_s=U_1B_a.
\]
Hence
\[
B_a(\ker(I-U_s))=\ker(I-U_1).
\]
By minimality of $\alpha$,
\[
\ker(I-U_1)=\Z 1_{\Omega_1}.
\]
Since $B_a^{-1}=B_a$, it follows that
\[
\ker(I-U_s)=\Z a\cong\Z.
\]

We now apply Proposition \ref{prop:coefficient-exact-sequences}. If $s$ is a sign coboundary, the
first short exact sequence becomes
\[
0\longrightarrow D_\alpha
\longrightarrow K_0(A)
\longrightarrow\Z
\longrightarrow0.
\]
Since $\Z$ is a free abelian group,
this sequence splits. Therefore
\[
K_0(A)\cong D_\alpha\oplus\Z.
\]
Moreover, by Proposition~ \ref{prop:coefficient-exact-sequences},
\[
K_1(A)\cong Q_s\oplus\Z.
\]
Since $Q_s$ is torsion-free by
Proposition \ref{prop:signed-torsion}, so is $K_1(A)$.

If $s$ is not a sign coboundary,
then $\ker(I-U_s)=0$, and the first
short exact sequence gives
\[
K_0(A)\cong D_\alpha.
\]
On the other hand,
Propositions \ref{prop:coefficient-exact-sequences} and \ref{prop:signed-torsion}
yield
\[\operatorname{Tor}(K_1(A))\cong\operatorname{Tor}(Q_s)\cong\Z/2\Z.\]
Since $(j_{\widetilde\alpha})_{*1}$
is injective on $Q_s$, the nonzero
torsion element is precisely
\[
(j_{\widetilde\alpha})_{*1}([v_s]).
\]
This completes the proof.
\end{proof}

\begin{rem}
We remind the readers that the isomorphisms involving $K_0(A)$ in Proposition \ref{prop:orientation-dichotomy}
are isomorphisms of abelian groups.
No assertion about the order structure
is made at this stage.
\end{rem}

\begin{cor}\label{cor:orientation-K1-invariant}
Let $\widetilde\alpha=\alpha\times\varphi$
and $\widetilde\beta=\beta\times\psi$
be minimal Cantor--circle extensions
on $X=\Omega_1\times\T$ and
$Y=\Omega_2\times\T$, respectively.
Let
\[
A=C(X)\rtimes_{\widetilde\alpha}\Z,
\qquad
B=C(Y)\rtimes_{\widetilde\beta}\Z.
\]
Suppose that
\[
K_1(A)\cong K_1(B)
\]
as abelian groups.

Then $s_{\widetilde\alpha}$ is a sign
coboundary over $\alpha$ if and only if
$s_{\widetilde\beta}$ is a sign
coboundary over $\beta$.

Equivalently, the continuous orientation
classes of $s_{\widetilde\alpha}$
and $s_{\widetilde\beta}$ are either
both trivial or both nontrivial.
\end{cor}

\begin{proof}
By Proposition~\ref{prop:orientation-dichotomy},
the orientation sign
$s_{\widetilde\alpha}$ is a sign
coboundary if and only if
$K_1(A)$ is torsion-free.
The same conclusion holds for
$s_{\widetilde\beta}$ and $K_1(B)$.

Since an isomorphism of abelian groups
preserves torsion subgroups, the
assumption $K_1(A)\cong K_1(B)$
implies that $K_1(A)$ is torsion-free
if and only if $K_1(B)$ is torsion-free.
The conclusion follows.
\end{proof}

\section{Base and orientation normalization}
\label{sec:normalization}

The finite coefficient $K$-condition in
Theorem~\ref{thm:main} prescribes the compatibility
of the induced $K$-theory maps associated with
the canonical embeddings of the coefficient
algebras into their crossed products.
By Proposition~\ref{prop:coefficient-exact-sequences},
the kernels of these maps are characterized by
ordinary and signed coboundaries. 

We begin by
formulating the resulting constraints on the
base homeomorphisms and fiber orientation signs in Lemma \ref{lem:coinvariant-compatibility}.

\subsection{Compatibility of coinvariant relations}
\label{subsec:coinvariant-compatibility}
Let $\taf=\alpha\times\varphi$ and
$\tbt=\beta\times\psi$ be minimal
Cantor--circle extensions, and write
\[
s=s_{\taf},\qquad t=s_{\tbt}.
\]
Suppose that a graded isomorphism
\[
k=(k_0,k_1):K_*(A)\longrightarrow K_*(B)
\]
and a sequence of unital $*$-isomorphisms
\[
\chi_n:C(X)\longrightarrow C(Y)
\]
satisfy the finite coefficient
$K$-condition
\eqref{eq:main-K-forward}.
Let
\[
\gamma_n:Y\longrightarrow X
\]
be the homeomorphism corresponding to
$\chi_n$, and write
\[
\gamma_n(\eta,z)
=
(g_n(\eta),q_{n,\eta}(z)),
\]
where each $g_n:\Omega_2\to\Omega_1$ is a homeomorphism. Put
\[
\varepsilon_n(\eta)
=
\deg(q_{n,\eta})\in\{-1,1\},\quad n\ge1.
\]

\begin{lem}\label{lem:coinvariant-compatibility}
For every pair of finite subsets
\[
F_0,F_1\subset C(\Omega_1,\Z),
\]
there exists $N\geq 1$ such that, for
all $n\geq N$, the following hold.

\begin{enumerate}
\item
For every $f\in F_0$,
\[
f\in(I-U_1^\alpha)C(\Omega_1,\Z)
\]
if and only if
\[
f\circ g_n
\in(I-U_1^\beta)C(\Omega_2,\Z).
\]

\item
For every $f\in F_1$,
\[
f\in(I-U_s)C(\Omega_1,\Z)
\]
if and only if
\[
\varepsilon_n\cdot(f\circ g_n)
\in(I-U_t)C(\Omega_2,\Z).
\]
\end{enumerate}
Here $U_1^\alpha$ and $U_1^\beta$
denote the canonical automorphisms
associated with $\alpha$ and $\beta$,
respectively.
\end{lem}

\begin{proof}
Choose finitely generated subgroups
\[
H_i\subset K_i(C(X))=C(\Omega_1,\Z),
\qquad i=0,1,
\]
containing the elements of $F_i$.

By \eqref{eq:main-K-forward}, there
exists $N\geq 1$ such that, for
$n\geq N$,
\[
k_i\big((j_{\taf})_{*i}(x)\big)
=
(j_{\tbt})_{*i}
\big((\chi_n)_{*i}(x)\big),
\qquad x\in H_i.
\]
Since $k_i$ is an isomorphism, we have
\[
(j_{\taf})_{*i}(x)=0
\quad\Longleftrightarrow\quad
(j_{\tbt})_{*i}
\big((\chi_n)_{*i}(x)\big)=0.
\]
Under the fiberwise rank identification,
\[
(\chi_n)_{*0}(f)=f\circ g_n.
\]
Moreover, Proposition~
\ref{prop:coefficient-exact-sequences}
gives
\[
\ker(j_{\taf})_{*0}
=
(I-U_1^\alpha)C(\Omega_1,\Z)
\]
and
\[
\ker(j_{\tbt})_{*0}
=
(I-U_1^\beta)C(\Omega_2,\Z).
\]
This proves the first assertion.

Similarly, we have
\[
(\chi_n)_{*1}(f)
=
\varepsilon_n\cdot(f\circ g_n).
\]
Applying Proposition~
\ref{prop:coefficient-exact-sequences}
to the two extensions gives
\[
\ker(j_{\taf})_{*1}
=
(I-U_s)C(\Omega_1,\Z)
\]
and
\[
\ker(j_{\tbt})_{*1}
=
(I-U_t)C(\Omega_2,\Z).
\]
The second assertion follows.
\end{proof}

\begin{rem}
Lemma~\ref{lem:coinvariant-compatibility} translates the finite
coefficient $K$-condition into finite-level constraints on the
underlying dynamics. By identifying the kernels of the
$K$-maps induced by coefficient algebras with ordinary and signed coboundaries, it shows that the
base maps $g_n$, together with their fiber-degree functions
$\varepsilon_n$, eventually preserve the relevant coboundary
relations on any prescribed finite set. Thus the $K_0$ condition
controls ordinary coinvariants, while the $K_1$ condition also
records the orientation twist. These are the compatibility
relations used in the subsequent normalization of the base dynamics
and fiber orientations; the lemma itself does not yet produce an
approximate conjugacy.
\end{rem}

\subsection{Normalization of the base dynamics}
\label{subsec:base-normalization}

We shall use the following finite-partition
matching result for Cantor minimal systems,
due to Lin and Matui
\cite[Lemma~4.7]{LinMatui2005AC};
see also \cite[Lemma~3.3]{Matui2005FullGroups}.

\begin{lem}[{\cite[Lemma~4.7]{LinMatui2005AC}}]
\label{lem:finite-partition-matching}
Let $\alpha$ and $\beta$ be minimal
homeomorphisms of a Cantor space $\Omega$,
and let $\mathcal P$ be a finite clopen
partition of $\Omega$.

Suppose that
\[
[1_U]=[1_{\beta(U)}]
\quad\text{in }D_\alpha
\]
for every $U\in\mathcal P$. Then there exists $\sigma\in[[\alpha]]$
such that
\[
\sigma\alpha\sigma^{-1}(U)=\beta(U),
\qquad U\in\mathcal P.
\]
\end{lem}

\begin{prop}\label{prop:base-normalization}
Under the assumptions of
Lemma~\ref{lem:coinvariant-compatibility}, let $\mathcal P$ be a finite clopen partition
of $\Omega_1$. Then, for all sufficiently
large $n$, there exists
$\sigma_n\in[[\beta]]$ such that
\[
g_n\sigma_n\beta\sigma_n^{-1} g_n^{-1}(U)
=
\alpha(U),
\qquad U\in\mathcal P.
\]
\end{prop}

\begin{proof}
For each $n$, put
\[
\delta_n=g_n^{-1}\alpha g_n,
\]
which is a minimal homeomorphism of
$\Omega_2$.

For every $U\in\mathcal P$, we have
\[
1_U-1_{\alpha(U)}
\in(I-U_1^\alpha)C(\Omega_1,\Z).
\]
By Lemma~\ref{lem:coinvariant-compatibility},
for all sufficiently large $n$,
\[
1_{g_n^{-1}(U)}
-
1_{g_n^{-1}(\alpha(U))}
\in(I-U_1^\beta)C(\Omega_2,\Z).
\]
Since
\[
g_n^{-1}(\alpha(U))
=
\delta_n(g_n^{-1}(U)),
\]
it follows that
\[
[1_V]=[1_{\delta_n(V)}]
\quad\text{in }D_\beta
\]
for every
$V\in g_n^{-1}(\mathcal P)$.

Applying Lemma~\ref{lem:finite-partition-matching}
to $\beta$, $\delta_n$, and
$g_n^{-1}(\mathcal P)$,
we obtain $\sigma_n\in[[\beta]]$ such that
\[
\sigma_n\beta\sigma_n^{-1}(V)
=
\delta_n(V)
\]
for every $V\in g_n^{-1}(\mathcal P)$. Conjugating by $g_n$ gives
\[
g_n\sigma_n\beta\sigma_n^{-1}g_n^{-1}(U)
=
\alpha(U)
\]
for every $U\in\mathcal P$.
\end{proof}

\begin{cor}\label{cor:base-approximate-conjugacy}
Suppose that the assumptions of
Theorem~\ref{thm:main}(2) hold. Then there exist sequences of homeomorphisms $h_n:\Omega_2\longrightarrow\Omega_1$ and 
$\ell_n:\Omega_1\longrightarrow\Omega_2$ such that
\[
d_\Homeo\bigl(h_n\beta h_n^{-1},\alpha\bigr)
\longrightarrow 0, \quad d_\Homeo\bigl(\ell_n\alpha\ell_n^{-1},\beta\bigr)
\longrightarrow 0.
\]
Moreover, the homeomorphisms $h_n$ may be
chosen in the form
\[
h_n=g_n\circ\sigma_n,
\qquad \sigma_n\in[[\beta]],
\]
where $g_n$ is associated with $\chi_n$.
\end{cor}

\begin{proof}
Choose a sequence of finite clopen
partitions $\{\mathcal P_m\}_{m\geq1}$
of $\Omega_1$ such that
\[
\operatorname{mesh}(\mathcal P_m)
\longrightarrow 0.
\]
By the uniform continuity of $\alpha$,
\[
\operatorname{mesh}(\alpha(\mathcal P_m))
\longrightarrow 0.
\]

For each $m$, Proposition~
\ref{prop:base-normalization} provides
an integer $N_m\geq1$ such that,
for every $n\geq N_m$, there exists
$\sigma_n\in[[\beta]]$ satisfying
\[
g_n\sigma_n\beta\sigma_n^{-1}g_n^{-1}(U)
=
\alpha(U),
\qquad U\in\mathcal P_m.
\]

Choose a strictly increasing sequence
of integers $n_m$ such that
$n_m\geq N_m$ for every $m$.
For $n_m\leq n<n_{m+1}$, choose
$\sigma_n\in[[\beta]]$ as above for
$\mathcal P_m$, and put
\[
h_n=g_n\circ\sigma_n,
\qquad
\delta_n=h_n\beta h_n^{-1}.
\]
The finitely many initial terms may
be chosen arbitrarily.

For every $U\in\mathcal P_m$, we have
\[
\delta_n(U)=\alpha(U).
\]
Consequently, for every $x\in U$,
\[
d(\delta_n(x),\alpha(x))
\leq\operatorname{diam}(\alpha(U)).
\]
Similarly, for every $y\in\alpha(U)$,
\[
d(\delta_n^{-1}(y),\alpha^{-1}(y))
\leq\operatorname{diam}(U).
\]
It follows that
\[
d_\Homeo(\delta_n,\alpha)
\leq
\operatorname{mesh}(\mathcal P_m)
+
\operatorname{mesh}(\alpha(\mathcal P_m)).
\]
The right-hand side tends to zero
as $n\to\infty$. Hence
\[
d_\Homeo(h_n\beta h_n^{-1},\alpha)
\longrightarrow0.
\]

Applying the same argument to the
sequence $\eta_n$ and the
equality \eqref{eq:main-K-backward}
gives homeomorphisms
$\ell_n:\Omega_1\to\Omega_2$ satisfying
\[
d_\Homeo(\ell_n\alpha\ell_n^{-1},\beta)
\longrightarrow0.
\]
The proof is complete.
\end{proof}

\begin{lem}\label{lem:full-group-inner-lift}
Let $\sigma\in[[\beta]]$, with
$\sigma(\eta)=\beta^{n(\eta)}(\eta)$. Then
\[
\Sigma(\eta,z)=\tbt^{n(\eta)}(\eta,z)
\]
is a homeomorphism covering $\sigma$. Put
$E_k=\{\eta:n(\eta)=k\}$ and
$p_k=j_{\tbt}(1_{E_k\times\T})$. The finite sum
\[
W_\sigma=\sum_k u_{\tbt}^kp_k
\]
is unitary and satisfies
\[
W_\sigma j_{\tbt}(f)W_\sigma^*=j_{\tbt}(f\circ\Sigma^{-1}).
\]
Consequently composition with $\Sigma$ or $\Sigma^{-1}$ preserves
the entire embedded coefficient $K$-maps and every coefficient trace
value. The same assertions hold for $\taf$.
\end{lem}
\begin{proof}
The sets $E_k$ and their images $\beta^k(E_k)$ form finite clopen
partitions. Iterates of $\tbt$ give the map and its inverse on these
pieces. 

Write $w_k=u_{\tbt}^kp_k$ and $q_k=u_{\tbt}^kp_ku_{\tbt}^{-k}$. Covariance gives
\[q_k=j_{\tbt}(1_{\bt^k(E_k)\times\T}),\quad w_k^*w_k=p_k,\quad w_kw_k^*=q_k.\]
Since $\sigma\in[[\bt]]$, one sees that the projections $p_k$ and $q_k$ each partition the identity. For $k\ne l$,
\[w_k^*w_l=w_k^*q_kq_lw_l=0,\qquad w_kw_l^*=w_kp_kp_lw_l^*=0.\]
Therefore,
\[W_\sigma W_\sigma^*=W_\sigma^* W_\sigma=1.\]
Moreover, 
\begin{align}\label{unitaryw}
    W_\sigma j_{\tbt}(f)W_\sigma^*&=\sum_k u_{\tbt}^kj_{\tbt}(1_{E_k\times\T}\cdot f)u_{\tbt}^{-k}\\
    &=j_{\tbt}\left(\sum_k (1_{E_k\times\T}\cdot f)\circ \tbt^{-k}\right).
\end{align}
On $\bt^k(E_k)\times\T$, the only nonzero summand is $f\circ \tbt^{-k}=f\circ\Sigma^{-1}$. This proves the equality. Finally, inner automorphisms fix
$K$-theory and traces, giving the remaining assertions.
\end{proof}

\begin{prop}\label{prop:combined-normalization}
Under the assumptions of Theorem~\ref{thm:main}\textup{(2)}, the
homeomorphisms associated with $\chi_n$ can be replaced by
$\widehat H_n:Y\to X$, with bases $h_n$, such that
\[
d_{\Homeo}(h_n\beta h_n^{-1},\alpha)\to0,\qquad
j_{\tbt}\widehat H_n^*=\Ad(W_n^*)j_{\tbt}\chi_n
\]
for unitaries $W_n\in B$.

In particular, the induced $K$-theory maps of the coefficient embeddings and the trace values on coefficient functions are preserved at every index. The same holds in the reverse direction.
\end{prop}

\begin{proof}
Use Corollary~\ref{cor:base-approximate-conjugacy} to write
$h_n=g_n\sigma_n$ with $\sigma_n\in[[\beta]]$.
Let $\Sigma_n$ be the corresponding extension from
Lemma~\ref{lem:full-group-inner-lift}, and let $H_n$ be the original
homeomorphism associated with $\chi_n$. Set $\widehat H_n=H_n\Sigma_n$.
Then
\[
j_{\tbt}\widehat H_n^*
=
\Ad(W_{\sigma_n}^*)j_{\tbt}\chi_n.
\]
Since inner automorphisms induce the identity on $K$-theory, for
$i=0,1$ we have
\[
(j_{\tbt}\widehat H_n^*)_{*i}
=
(j_{\tbt}\chi_n)_{*i}.
\]
Moreover, for every $f\in C(\Omega_1\times\T)$ and $\tau\in T(B)$,
\[
\tau(j_{\tbt}\widehat H_n^*(f))
=
\tau(j_{\tbt}\chi_n(f)),
\]
because tracial states are invariant under inner automorphisms.
These equalities hold for each $n$, without passing to a limit. The reverse direction follows by applying the same construction to
$\eta_n$.
\end{proof}

\subsection{Orbit-breaking signed coinvariants}
\label{subsec:broken-signed-coinvariants}
The normalization established in
Proposition~\ref{prop:combined-normalization}
aligns the base dynamics while preserving
the coefficient-induced $K$-theory maps.
It does not, however, ensure that the
orientation signs of the conjugated
skew products agree pointwise with those
of the target systems.

To obtain exact orientation matching,
we first consider the case of a nontrivial
continuous orientation class.
The key additional ingredient is an
orbit-breaking version of the signed
coinvariant group, in which the coboundary
relations are restricted at a fixed base point.
This refinement records the integer ambiguity
that remains when relations in the ordinary
signed coinvariant group are lifted to the
orbit-breaking quotient.
We shall use local corrections from the
topological full group to eliminate this
ambiguity.

Let $(\Omega,\beta)$ be a Cantor minimal
system, and let
$t:\Omega\to\{-1,1\}$ be continuous.
Fix $\eta_0\in\Omega$. Put $M=C(\Omega,\mathbb Z)$ as before and
\[
M_{\eta_0}
=
\{f\in M:f(\eta_0)=0\}.
\]
Recall the automorphism
\[
(U_t f)(\eta)
=
t(\beta^{-1}(\eta))f(\beta^{-1}(\eta))
\]
for $\eta\in\Omega$. We define the orbit-breaking signed coinvariant
group at $\eta_0\in\Omega$ by
\[
Q_{t,\eta_0}^{\mathrm{br}}
=
M/(I-U_t)M_{\eta_0}.
\]
There is a canonical surjective
homomorphism
\[
\pi_{\eta_0}:
Q_{t,\eta_0}^{\mathrm{br}}
\longrightarrow Q_t.
\]
When necessary, we write
$[f]_{\mathrm{br}}$ for the class of
$f\in M$ in $Q_{t,\eta_0}^{\mathrm{br}}$,
and $[f]_t$ for its class in $Q_t$.
For convenience, we shall also use the integer-valued
function
\[
h_t
=
\frac{1_\Omega-t\circ\beta^{-1}}{2}
\in C(\Omega,\mathbb Z),
\]
so that
\[
(I-U_t)1_\Omega=2h_t.
\]

\begin{lem}\label{lem:broken-signed-kernel}
Suppose that $t$ is not a continuous
sign coboundary over $\beta$.
Then the kernel of the canonical map
\[
\pi_{\eta_0}:
Q_{t,\eta_0}^{\mathrm{br}}
\longrightarrow Q_t
\]
is infinite cyclic, generated by
\[[2h_t]_{\mathrm{br}}.\]
\end{lem}

\begin{proof}
We first show that $I-U_t$ is injective.
Suppose $U_t f=f$ for some
$f\in C(\Omega,\mathbb Z)$.
Then $|f|$ is $\beta$-invariant,
and hence constant by minimality.
If $f\ne 0$, the function
$f/|f|$ is a continuous sign
satisfying
\[
\frac{f(\beta\eta)}{|f(\beta\eta)|}
=
t(\eta)\frac{f(\eta)}{|f(\eta)|},
\]
contrary to the assumption.
Thus $\ker(I-U_t)=0$.

Every $g\in M$ can be written as
\[
g=g(\eta_0)1_\Omega+
\bigl(g-g(\eta_0)1_\Omega\bigr),
\]
where the second summand belongs
to $M_{\eta_0}$.
Consequently, the kernel of
$\pi_{\eta_0}$ is generated by
$[(I-U_t)1_\Omega]$.
Finally, suppose
\[
k(I-U_t)1_\Omega=(I-U_t)g
\]
for some $g\in M_{\eta_0}$
and $k\in\mathbb Z$.
Injectivity of $I-U_t$ gives
$g=k1_\Omega$.
Evaluating at $\eta_0$ yields
$k=0$.
Thus the generator has infinite order.
\end{proof}

\medskip

Now we take a Kakutani--Rokhlin partition
\[
\mathcal P
=
\{V_{i,j}:
1\le i\le m,\ 1\le j\le r_i\},
\]
where $V_{i,j}=\beta^{j-1}(V_i)$ with $V_i$ being the bases, and suppose that the union of its
roof levels is a clopen neighborhood
$R$ of $\eta_0$.

\begin{lem}\label{lem:broken-signed-tower-sums}
Let $f\in C(\Omega,\mathbb Z)$ and $t:\Omega\to\{-1,1\}$ be continuous. Assume that $f$ and $t$ are
constant on every level $V_{i,j}$.
Write
\[
f_{i,j}=f|_{V_{i,j}},
\qquad
t_{i,j}=t|_{V_{i,j}},
\]
and define
\[
b_{i,j}
=
\prod_{\ell=j}^{r_i-1}t_{i,\ell},
\qquad
b_{i,r_i}=1.
\]
For each tower, put
\[
T_i(f)
=
\sum_{j=1}^{r_i}b_{i,j}f_{i,j}.
\]
Then
\begin{equation}
\label{eq:broken-signed-tower-reduction}
[f]_{\mathrm{br}}
=
\left[
\sum_{i=1}^m
T_i(f)1_{V_{i,r_i}}
\right]_{\mathrm{br}}.
\end{equation}
Moreover, the following are equivalent:
\begin{enumerate}
\item
$[f]_{\mathrm{br}}=0$;

\item
there exists a Kakutani--Rokhlin
partition as above, with roof union
containing $\eta_0$, such that
\[
T_i(f)=0,
\qquad 1\le i\le m.
\]
\end{enumerate}
\end{lem}

\begin{proof}
We first establish
\eqref{eq:broken-signed-tower-reduction}.

For $1\le i\le m$ and
$1\le j<r_i$, the function
$1_{V_{i,j}}$ vanishes at $\eta_0$,
since $\eta_0$ belongs to the
union of the roof levels.
Consequently,
\[
[(I-U_t)1_{V_{i,j}}]_{\mathrm{br}}=0.
\]
Since $t$ is constant on $V_{i,j}$,
we have
\[
U_t1_{V_{i,j}}
=
t_{i,j}1_{V_{i,j+1}}.
\]
It follows that
\[
[1_{V_{i,j}}]_{\mathrm{br}}
=
t_{i,j}
[1_{V_{i,j+1}}]_{\mathrm{br}}.
\]
Iterating this identity gives
\[
[1_{V_{i,j}}]_{\mathrm{br}}
=
b_{i,j}
[1_{V_{i,r_i}}]_{\mathrm{br}}.
\]
Since
\[
f=
\sum_{i=1}^m
\sum_{j=1}^{r_i}
f_{i,j}1_{V_{i,j}},
\]
we obtain
\eqref{eq:broken-signed-tower-reduction}. We next prove the equivalence.

Suppose first that
$[f]_{\mathrm{br}}=0$.
Then there exists
$g\in C(\Omega,\mathbb Z)$
with $g(\eta_0)=0$ such that
\[
f=(I-U_t)g.
\]
Since $g$ is locally constant,
there is a clopen neighborhood
$R_0$ of $\eta_0$ on which
$g$ vanishes. Choose a clopen neighborhood
$R\subseteq R_0$ of $\eta_0$,
and form the first-return
Kakutani--Rokhlin partition
with base $\beta(R)$.
Refining its bases if necessary,
we may assume that $f$, $g$,
and $t$ are constant on every
tower level. Write
\[
g_{i,j}=g|_{V_{i,j}}.
\]
Since $g$ vanishes on the
roof union $R$, we have
\[
g_{i,r_i}=0.
\]
The identity $f=(I-U_t)g$
gives
\[
f_{i,1}=g_{i,1}
\]
and
\[
f_{i,j}
=
g_{i,j}
-
t_{i,j-1}g_{i,j-1},
\qquad 2\le j\le r_i.
\]
Using
\[
b_{i,j-1}
=
t_{i,j-1}b_{i,j},
\]
the signed sum telescopes:
\[
\begin{aligned}
T_i(f)
&=
\sum_{j=1}^{r_i}
b_{i,j}f_{i,j}\\
&=
g_{i,r_i}\\
&=0.
\end{aligned}
\]
Conversely, suppose that a
partition as in the statement
satisfies
\[
T_i(f)=0
\qquad(1\le i\le m).
\]
Define integers $g_{i,j}$
recursively by
\[
g_{i,1}=f_{i,1},
\]
and
\[
g_{i,j}
=
f_{i,j}
+
t_{i,j-1}g_{i,j-1},
\qquad 2\le j\le r_i.
\]
By induction,
\[
g_{i,r_i}=T_i(f)=0.
\]
Define $g\in C(\Omega,\mathbb Z)$
by setting
\[
g|_{V_{i,j}}=g_{i,j}.
\]
This is continuous because the
tower levels form a finite
clopen partition.
Moreover, $g$ vanishes on
the roof union $R$, and hence
$g(\eta_0)=0$. The defining recurrence shows that
\[
f=(I-U_t)g
\]
at every non-base level.
At a base level, its predecessor
belongs to $R$, where $g=0$,
so the same identity holds there.
As a conclusion, $g\in M_{\eta_0}$ and
$[f]_{\mathrm{br}}=0$.

Finally, the same argument applies simultaneously to
any finite collection of such equalities.
Indeed, one may choose a common clopen
roof neighborhood on which all the
corresponding coboundary witnesses vanish,
and then take a common finite refinement
of the associated tower partitions.
\end{proof}

The following lemma provides a local correction mechanism
for the signed coinvariant obstruction. It realizes
arbitrary integer multiples of $[h_t]_{\mathrm{br}}$
by characteristic functions supported in any prescribed
nonempty clopen set, together with a compatible signed
pushforward induced by an element of the topological
full group.

\begin{lem}\label{lem:local-signed-correction}
Assume that $t$ is not a sign coboundary. For every integer $m$ and every
nonempty clopen set $O\subset\Omega$, there are a clopen $E\subset O$ and
$\rho\in[[\beta]]$, equal to the identity off $E$, whose signed pushforward
\[
S_\rho f=\sum_k U_t^k(f1_{\{n_\rho=k\}})
\]
satisfies
\[
[1_E]_{\mathrm{br}}=m[h_t]_{\mathrm{br}},\qquad S_\rho1_E=-1_E.
\]
\end{lem}

\begin{proof}
We first prove the assertion for
$m=\epsilon\in\{-1,1\}$. By Lemma~\ref{lem:orientation-cover}, the sign extension
\[
D(\eta,e)=(\beta(\eta),t(\eta)e),
\qquad e\in\{-1,1\},
\]
is minimal.

\smallskip
\noindent
\textbf{Step 1. A uniform backward hitting-time bound.}
Since $D$ is minimal and
$O\times\{\epsilon\}$ is a nonempty open subset
of $\Omega\times\{-1,1\}$, we have
\[
\Omega\times\{-1,1\}
=
\bigcup_{k\geq 0}
D^k(O\times\{\epsilon\}).
\]
By compactness, there exists $N\geq 1$ such that
\begin{equation}\label{eq:local-signed-hitting}
\Omega\times\{-1,1\}
=
\bigcup_{k=0}^{N}
D^k(O\times\{\epsilon\}).
\end{equation}
In particular, for every $\eta\in\Omega$,
there exists $0\leq k\leq N$ such that
\[
D^{-k}(\eta,1)\in O\times\{\epsilon\}.
\]
Choose a sufficiently small clopen neighborhood
$R$ of $\eta_0$ such that $t$ is constant
on $R$, with value $c\in\{-1,1\}$, and
\[
R\cap\beta^k(R)=\varnothing,
\qquad 1\leq k\leq N.
\]
By minimality, the 
Kakutani--Rokhlin partition constructed from the first return time with roof
union $R$ has all tower heights greater
than $N$. 
Denote the resulting partition by
\[
\mathcal P
=
\{V_{i,j}:1\le i\le q,\,
1\le j\le r_i\}.
\]
Refining its bases if necessary,
we may assume that $t$ is constant on
every tower level.
Write
\[
t_{i,j}=t|_{V_{i,j}},
\qquad
b_{i,j}
=
\prod_{\ell=j}^{r_i-1}t_{i,\ell},
\qquad
b_{i,r_i}=1,
\]
and put
\[
T_i=\prod_{j=1}^{r_i}t_{i,j}
=cb_{i,1}.
\]

\smallskip
\noindent
\textbf{Step 2. Computation of the signed
tower sums.} Recall that
\[
h_t=\frac{1-t\circ\beta^{-1}}2.
\]
Since the predecessor of each base level
belongs to $R$, we have
\[
h_t|_{V_{i,1}}=\frac{1-c}{2}.
\]
For $2\leq j\leq r_i$,
\[
h_t|_{V_{i,j}}
=
\frac{1-t_{i,j-1}}2.
\]
By definition,
\[
T_i(h_t)
=
\sum_{j=1}^{r_i}
b_{i,j}h_t|_{V_{i,j}}.
\]
Upon using
\[
b_{i,j-1}=t_{i,j-1}b_{i,j},
\]
we obtain
\[
\begin{aligned}
2T_i(h_t)
&=
b_{i,1}(1-c)
+
\sum_{j=2}^{r_i}
b_{i,j}(1-t_{i,j-1})\\
&=
b_{i,1}(1-c)
+
\sum_{j=2}^{r_i}
(b_{i,j}-b_{i,j-1})\\
&=
b_{i,1}(1-c)
+b_{i,r_i}-b_{i,1}\\
&=
1-cb_{i,1}\\
&=1-T_i.
\end{aligned}
\]
Consequently,
\begin{equation}\label{eq:torsion-tower-sum}
T_i(h_t)=\frac{1-T_i}{2},
\qquad
[h_t]_{\mathrm{br}}
=
\left[
\sum_{T_i=-1}1_{V_{i,r_i}}
\right]_{\mathrm{br}}.
\end{equation}
The second identity follows from the
signed tower representation in
Lemma~\ref{lem:broken-signed-tower-sums}.

\smallskip
\noindent
\textbf{Step 3. Construction of the
correction set.} Fix a tower with $T_i=-1$.
For each point $\eta$ in its roof,
\eqref{eq:local-signed-hitting}
provides $0\leq k\leq N$ such that
\[
D^{-k}(\eta,1)
\in O\times\{\epsilon\}.
\]
Since the backward segment starts at a
roof level and $r_i>N$, all of its
points lie in the same tower.

If its endpoint belongs to level
$V_{i,j}$, then the product of the
orientation signs from that level
to the roof is $b_{i,j}$.
The second coordinate of the backward
iterate is therefore $b_{i,j}$.

Thus, for every point in the roof of
an odd-sign tower, there is a backward
orbit segment of length at most $N$
whose endpoint belongs to $O$ and whose
orientation sign is $\epsilon$.

Since $O$ is clopen and only finitely
many hitting times are involved, we may
refine the tower bases into finitely many
clopen subsets on each of which a
suitable hitting time is constant.
Consequently, on every refined tower
with $T_i=-1$, there is a selected
level $V_{i,k_i}$ satisfying
\[
V_{i,k_i}\subset O,
\qquad
b_{i,k_i}=\epsilon.
\]

There are only finitely many possible
hitting times. By partitioning the roof
according to the first admissible
hitting time and pulling these clopen
sets back to the tower bases, we obtain
a finite clopen refinement for which
a selected index $k_i$ is constant
on each refined tower and satisfies
\[
V_{i,k_i}\subset O,
\qquad
b_{i,k_i}=\epsilon.
\]
We make this selection in every tower
with $T_i=-1$, and make no selection
in towers with $T_i=1$. Now we let
\[
E_\epsilon
=
\bigsqcup_{T_i=-1}V_{i,k_i}.
\]
Then $E_\epsilon$ is clopen and
$E_\epsilon\subset O$. For each $i$,
\[
T_i(1_{E_\epsilon})
=
\begin{cases}
\epsilon,&T_i=-1,\\
0,&T_i=1.
\end{cases}
\]
Comparing this with
\eqref{eq:torsion-tower-sum}, we obtain
\[
T_i(1_{E_\epsilon})
=
\epsilon T_i(h_t)
\qquad(1\leq i\leq q).
\]
The signed tower criterion in
Lemma~\ref{lem:broken-signed-tower-sums}
therefore gives
\begin{equation}\label{eq:local-signed-class}
[1_{E_\epsilon}]_{\mathrm{br}}
=
\epsilon[h_t]_{\mathrm{br}}.
\end{equation}
Moreover, at least one tower has
$T_i=-1$. Otherwise
\eqref{eq:torsion-tower-sum} would imply
$[h_t]_{\mathrm{br}}=0$, contrary to
Lemma~\ref{lem:broken-signed-kernel}.
In particular, $E_\epsilon$ is nonempty.

\smallskip
\noindent
\textbf{Step 4. Construction of the
full-group element.} Let $r_{E_\epsilon}(\eta)$ be the first return time of $\eta\in E_\epsilon$
to $E_\epsilon$. Since $\beta$ is minimal
and $E_\epsilon$ is nonempty and clopen,
$r_{E_\epsilon}$ is bounded
and locally constant. Define
\[
\rho_\epsilon(\eta)
=
\begin{cases}
\beta^{r_{E_\epsilon}(\eta)}(\eta),
&\eta\in E_\epsilon,\\
\eta,&\eta\notin E_\epsilon.
\end{cases}
\]
The first-return map is a homeomorphism
of $E_\epsilon$, so its extension by
the identity is an element of
$[[\beta]]$. In particular,
\[
\rho_\epsilon(E_\epsilon)=E_\epsilon,
\qquad
\rho_\epsilon|_{\Omega\setminus E_\epsilon}
=\id.
\]

We claim that
\begin{equation}\label{eq:local-signed-return}
\prod_{\ell=0}^{r_{E_\epsilon}(\eta)-1}
t(\beta^\ell\eta)=-1
\qquad(\eta\in E_\epsilon).
\end{equation}
In fact, let $\eta$ belong to the
selected level of an odd-sign tower.
The product of the signs from $\eta$
to the roof of that tower, excluding
the roof sign itself, is
\[
b_{i,k_i}=\epsilon.
\]
Before the orbit reaches the next selected level, it may pass through
several complete towers. Every such
 tower must have total sign
$1$, since a tower of total sign $-1$
contains a selected level, which would
give an earlier return to $E_\epsilon$.

Suppose that the next selected level
belongs to a tower indexed by $i'$.
The product of the signs from the
preceding roof, including its sign,
to that selected level, excluding the
sign at the selected level, is
\[
c\prod_{j=1}^{k_{i'}-1}t_{i',j}
=
\frac{T_{i'}}{b_{i',k_{i'}}}
=
-\epsilon.
\]
Thus the sign of the complete
first-return segment is
\[
\epsilon\cdot 1\cdot(-\epsilon)=-1.
\]
This proves
\eqref{eq:local-signed-return}.

By the definition of signed pushforward,
the value of
$S_{\rho_\epsilon}1_{E_\epsilon}$
at a point of $E_\epsilon$ is the
orientation sign of the unique
return segment ending at that point.
Equation~\eqref{eq:local-signed-return}
therefore implies
\begin{equation}\label{eq:local-signed-pushforward}
S_{\rho_\epsilon}1_{E_\epsilon}
=
-1_{E_\epsilon}.
\end{equation}
Together,
\eqref{eq:local-signed-class} and
\eqref{eq:local-signed-pushforward}
prove the assertion for
$m=\epsilon\in\{-1,1\}$.

\smallskip
\noindent
\textbf{Step 5. Arbitrary integer multiples.} If $m=0$, take
\[
E=\varnothing,\qquad \rho=\id.
\]
Suppose that $m\neq0$, and put
\[
p=|m|,
\qquad
\epsilon=\operatorname{sign}(m).
\]
Since $\Omega$ is a Cantor space
and $O$ is nonempty and clopen,
we may choose pairwise disjoint
nonempty clopen sets
\[
O_1,\ldots,O_p\subset O.
\]
Applying Steps 1--4 separately to
each $O_j$, we obtain clopen sets
$E_j\subset O_j$ and elements
$\rho_j\in[[\beta]]$ such that
\[
[1_{E_j}]_{\mathrm{br}}
=
\epsilon[h_t]_{\mathrm{br}},
\qquad
S_{\rho_j}1_{E_j}=-1_{E_j},
\]
and
\[
\rho_j(E_j)=E_j,
\qquad
\rho_j|_{\Omega\setminus E_j}=\id.
\]
Set
\[
E=\bigsqcup_{j=1}^{p}E_j,
\qquad
\rho=\rho_1\circ\cdots\circ\rho_p.
\]
Because the $E_j$ are pairwise
disjoint and invariant under their
respective maps, the $\rho_j$ have
pairwise disjoint invariant supports.
They therefore commute, and $\rho$
belongs to $[[\beta]]$.
Moreover,
\[
\rho|_{\Omega\setminus E}=\id,
\qquad
\rho|_{E_j}=\rho_j|_{E_j}.
\]
Additivity in the broken signed
coinvariant group gives
\[
\begin{aligned}
[1_E]_{\mathrm{br}}
&=
\sum_{j=1}^{p}[1_{E_j}]_{\mathrm{br}}\\
&=
p\epsilon[h_t]_{\mathrm{br}}\\
&=
m[h_t]_{\mathrm{br}}.
\end{aligned}
\]
Finally, the signed pushforward of
$\rho$ restricts on each $E_j$ to
that of $\rho_j$. Hence
\[
\begin{aligned}
S_\rho1_E
&=
\sum_{j=1}^{p}S_\rho1_{E_j}\\
&=
\sum_{j=1}^{p}S_{\rho_j}1_{E_j}\\
&=
-\sum_{j=1}^{p}1_{E_j}\\
&=-1_E.
\end{aligned}
\]
This completes the proof.
\end{proof}

\subsection{Simultaneous normalization of base and orientation}
\label{subsec:simultaneous-normalization}
We first establish the finite signed alignment when the orientation classes
are nontrivial.

\begin{prop}\label{prop:finite-signed-alignment}
Suppose $s=s_{\taf}$ and $t=s_{\tbt}$ are not sign coboundaries.
Let $\lambda:K_*(B)\to K_*(A)$ be a graded isomorphism and let
$\xi:C(Y)\to C(X)$ be an isomorphism. Choose a first-return partition
\[
\mathcal Q
=
\{W_{w,l}:1\le w\le q,\ 1\le l\le h_w\}
\]
for $\beta$, with $t$ constant on
its levels and with one common value $\tau_0$ on the roof union $R$.
Let $\mathcal Q^\flat$ consist of the nonroof levels and the merged roof $R$.
Assume, in both degrees, that
\[
(j_{\taf}\xi)_{*i}(1_U)=\lambda_i(j_{\tbt})_{*i}(1_U)
\quad(U\in\mathcal Q,\ i=0,1).
\]
Then there exist $\sigma\in[[\alpha]]$
and its iterate lift $\Sigma: \Omega_1\times\T\to \Omega_1\times\T$
such that, upon setting
\[
\xi'(f)=\xi(f)\circ\Sigma^{-1},
\]
the associated homeomorphism $H':X\to Y$
and its base map $h':\Omega_1\to\Omega_2$
satisfy
\[
h'\alpha(h')^{-1}(U)=\beta(U)\quad(U\in\mathcal Q^\flat),\qquad
\deg(H'\taf(H')^{-1})_\eta=t(\eta).
\]
Moreover $j_{\taf}\xi'=\Ad(W)j_{\taf}\xi$ for a unitary $W\in A$.
\end{prop}

\begin{proof}
Let
\[
H(\omega,z)
=
(\gamma^{-1}(\omega),q_\omega(z))
\]
be the homeomorphism associated with $\xi$, and put
\[
a(\omega)=\deg(q_\omega)\in\{-1,1\}.
\]
Denote the homomorphism between $K$-theory induced from $\xi$ by
\[
L_0f=f\circ\gamma^{-1},
\qquad
L_1f=a(f\circ\gamma^{-1}).
\]
We now fix $\omega_0\in\Omega_1$, and use the orbit-breaking signed
coinvariant group associated with $(\alpha,s,\omega_0)$. We first modify $H$ by iterate corrections
so that the required relations hold in the orbit-breaking
signed quotient. We then realize these relations
by a finite permutation of return-tower levels.

\smallskip
\noindent
\textbf{Step 1. Alignment in the orbit-breaking signed quotient.} Write
\[
t_{w,l}=t|_{W_{w,l}},
\qquad
T_w=\prod_{l=1}^{h_w}t_{w,l},
\]
and define
\[
g=\sum_{T_w=-1}1_{W_{w,h_w}}.
\]
By the signed tower computation
\eqref{eq:torsion-tower-sum}, applied to
$(\beta,t)$, we have
\[
[g]_t=[h_t]_t,
\qquad
h_t=\frac{1-t\circ\beta^{-1}}2.
\]
Since $t$ is not a sign coboundary,
Proposition~\ref{prop:signed-torsion}
shows that this is the unique nonzero
torsion class in the signed coinvariant group.
The compatibility hypothesis gives
\[
[L_1g]_s=[h_s]_s,
\qquad
h_s=\frac{1-s\circ\alpha^{-1}}2.
\]
Here the equality follows by applying the
hypothesis to the finitely many roof indices
occurring in $g$.

By Lemma~\ref{lem:broken-signed-kernel},
the kernel of the natural homomorphism
from the orbit-breaking signed quotient to the
ordinary signed quotient is generated
by $[2h_s]_{\mathrm{br}}$. Consequently,
there is an integer $n$ such that
\begin{equation}\label{eq:finite-alignment-initial}
[L_1g]_{\mathrm{br}}
=
(2n+1)[h_s]_{\mathrm{br}}.
\end{equation}
There is at least one tower with $T_w=-1$.
Choose such a tower and a nonempty clopen set
\[
O\subset\gamma(W_{w,h_w})
\]
on which $a$ is constant, say $a|_O=\epsilon$,
where $\epsilon\in\{-1,1\}$.

Apply Lemma~\ref{lem:local-signed-correction}
with $m=n\epsilon$. We obtain a clopen
set $E\subset O$ and $\rho\in[[\alpha]]$
such that
\[
[1_E]_{\mathrm{br}}
=
n\epsilon[h_s]_{\mathrm{br}},
\qquad
S_\rho1_E=-1_E,
\]
with $\rho$ equal to the identity outside $E$. Let $\Sigma_\rho$ be the iterate lift
of $\rho$, and replace $\xi$ by
\[
\xi_\rho(f)=\xi(f)\circ\Sigma_\rho^{-1}.
\]
The new winding map is $S_\rho L_1$.
Since $E$ is contained in the image of
a single roof label, and
$L_1g=\epsilon$ on $E$, we have
\[
S_\rho L_1g
=
L_1g-\epsilon1_E
+\epsilon S_\rho1_E
=
L_1g-2\epsilon1_E.
\]
It follows from
\eqref{eq:finite-alignment-initial} that
\[
\begin{aligned}
[S_\rho L_1g]_{\mathrm{br}}
&=
(2n+1)[h_s]_{\mathrm{br}}
-2\epsilon[1_E]_{\mathrm{br}}\\
&=
(2n+1)[h_s]_{\mathrm{br}}
-2n[h_s]_{\mathrm{br}}\\
&=[h_s]_{\mathrm{br}}.
\end{aligned}
\]
Moreover, $\rho$ preserves the image of
every source label set. Thus every unsigned
label function is unchanged, and the signed
label functions outside the selected roof
label are unchanged.

We next impose the signed relations between
consecutive levels of each source tower.
Fix $w$ and proceed in descending order
\[
l=h_w-1,h_w-2,\ldots,1.
\]
For $U=W_{w,l}$, covariance gives
\[
U_t1_U=t_{w,l}1_{W_{w,l+1}}.
\]
The compatibility hypothesis therefore implies
\[
[L_11_U]_s
=
t_{w,l}[L_11_{W_{w,l+1}}]_s.
\]
The same equality remains valid after the
preceding full-group corrections, since their
signed pushforwards induce the identity on
the ordinary signed coinvariant group.
Hence, for some $d\in\mathbb Z$,
\[
[L_11_U]_{\mathrm{br}}
-
t_{w,l}[L_11_{W_{w,l+1}}]_{\mathrm{br}}
=
2d[h_s]_{\mathrm{br}}.
\]

Choose a nonempty clopen subset of the 
image of $U$ on which the fiber degree
is a constant $\epsilon\in\{-1,1\}$.
Applying
Lemma~\ref{lem:local-signed-correction}
with $m=d\epsilon$ gives a correction
supported in that image, changing
$L_11_U$ by $-2\epsilon1_E$. Its orbit-breaking class therefore changes by
\[
-2\epsilon[1_E]_{\mathrm{br}}
=
-2d[h_s]_{\mathrm{br}},
\]
and the desired equality holds.

Since the correction preserves the image
of $U$ and fixes every other label
function, it does not alter the equalities
already established at higher levels.
In particular, corrections supported in
nonroof labels leave $L_1g$ unchanged.
After finitely many such corrections,
we obtain
\begin{equation}\label{eq:broken-alignment-relations}
\begin{aligned}
[L_1g]_{\mathrm{br}}
&=[h_s]_{\mathrm{br}},\\
[L_11_{W_{w,l}}]_{\mathrm{br}}
&=
t_{w,l}
[L_11_{W_{w,l+1}}]_{\mathrm{br}},
\qquad 1\leq l<h_w.
\end{aligned}
\end{equation}
Let $\rho_{\mathrm{tot}}$ denote the
composition of all corrections made so far.
The current base map from $\Omega_2$
to $\Omega_1$ is
\[
\gamma_c=\rho_{\mathrm{tot}}\circ\gamma.
\]
Each correction preserves every
label set, so
\[
\gamma_c(U)=\gamma(U)
\qquad(U\in\mathcal Q).
\]
From now on, $\gamma$, $a$, $L_0$, and $L_1$
denote these corrected objects. Moreover, by Lemma~\ref{lem:full-group-inner-lift},
all corrections in this step are implemented
by inner automorphisms of $A$.

\smallskip
\noindent
\textbf{Step 2. Exact tower-sum identities.} For every $w$ and $l<h_w$, the degree-zero
compatibility hypothesis gives
\[
[L_01_{W_{w,l}}]_\alpha
=
[L_01_{W_{w,l+1}}]_\alpha.
\]
In particular, the full group corrections of Step 1
preserve their ordinary coinvariant classes. Consequently, there are continuous
integer-valued functions $u_{w,l}$ such that
\begin{equation}\label{eq:alignment-unsigned-witness}
L_01_{W_{w,l}}
-
L_01_{W_{w,l+1}}
=
(I-U_1)u_{w,l}.
\end{equation}
Subtracting the constant $u_{w,l}(\omega_0)$
does not change this identity.
We may therefore assume
\[
u_{w,l}(\omega_0)=0.
\]
The orbit-breaking signed equalities
\eqref{eq:broken-alignment-relations}
also provide continuous integer-valued
functions $v_{w,l}$ and $v_0$,
all vanishing at $\omega_0$, such that
\begin{equation}\label{eq:alignment-signed-witness}
L_11_{W_{w,l}}
-
t_{w,l}L_11_{W_{w,l+1}}
=
(I-U_s)v_{w,l}
\end{equation}
and
\begin{equation}\label{eq:alignment-roof-witness}
L_1g-h_s=(I-U_s)v_0.
\end{equation}
Choose a sufficiently small clopen
neighborhood $R_\alpha$ of $\omega_0$
on which $s$ is constant, with value
$c\in\{-1,1\}$, and on which all the
functions $u_{w,l}$, $v_{w,l}$, and
$v_0$ vanish.
We may additionally require that the
first-return towers associated with
$R_\alpha$ have arbitrarily large heights.
Refine their bases so that $s$, $a$,
and every pulled-back source label
are constant on each level. Denote the resulting partition by
\[
\mathcal P
=
\{Z_{v,j}:1\leq v\leq p,\,
1\leq j\leq r_v\}.
\]
Its roof union is $R_\alpha$. To avoid ambiguity, we call the towers in $\mathcal P$ the auxiliary towers.

Choose a source tower $w_*$ with
$T_{w_*}=-1$. Since
$\gamma(W_{w_*,h_{w_*}})$ is nonempty
and clopen, minimality of $\alpha$
gives a uniform bound on forward
hitting times to this set.
Taking the auxiliary tower heights
sufficiently large ensures that
every auxiliary tower contains at
least one level labeled
$(w_*,h_{w_*})$. For each auxiliary tower, write
\[
s_{v,j}=s|_{Z_{v,j}},
\qquad
b_{v,j}
=
\prod_{k=j}^{r_v-1}s_{v,k},
\qquad
b_{v,r_v}=1,
\]
and put
\[
S_v=cb_{v,1}.
\]
Define the signed color of its
$j$-th level by
\[
e_{v,j}
=
a|_{Z_{v,j}}b_{v,j}
\in\{-1,1\}.
\]
For a function constant on the levels
of an auxiliary tower, write
\[
\mathcal T_v^s(f)
=
\sum_{j=1}^{r_v}b_{v,j}f|_{Z_{v,j}},
\qquad
\mathcal T_v^1(f)
=
\sum_{j=1}^{r_v}f|_{Z_{v,j}}.
\]
Since the witnesses in
\eqref{eq:alignment-unsigned-witness}--
\eqref{eq:alignment-roof-witness}
vanish on the roof union,
telescoping along each tower gives
\begin{equation}\label{eq:alignment-exact-sums}
\begin{aligned}
\mathcal T_v^1(L_01_{W_{w,l}})
&=
\mathcal T_v^1(L_01_{W_{w,l+1}}),\\
\mathcal T_v^s(L_11_{W_{w,l}})
&=
t_{w,l}
\mathcal T_v^s(L_11_{W_{w,l+1}}),\\
\mathcal T_v^s(L_1g)
&=
\mathcal T_v^s(h_s)
=
\frac{1-S_v}{2}.
\end{aligned}
\end{equation}
These are exact integer equalities
on every auxiliary tower.

\smallskip
\noindent
\textbf{Step 3. Matching the colored levels.} 
Each level $Z_{v,j}$ of the auxiliary partition
is contained in a unique set $\gamma(W_{w,l})$.
We assign to $Z_{v,j}$ the label $(w,l)$
determined by this inclusion.
For $e\in\{-1,1\}$, define
\[
N(v,w,l,e)
=
\#\{j:
Z_{v,j}\subset\gamma(W_{w,l}),
\ e_{v,j}=e\}.
\]
The first equality in
\eqref{eq:alignment-exact-sums}
gives
\begin{equation}\label{eq:alignment-count-total}
\sum_{e\in\{-1,1\}}N(v,w,l,e)
=
\sum_{e\in\{-1,1\}}N(v,w,l+1,e).
\end{equation}
Since
\[
L_11_{W_{w,l}}
=
a1_{\gamma(W_{w,l})},
\]
the second equality in
\eqref{eq:alignment-exact-sums}
gives
\begin{equation}\label{eq:alignment-count-signed}
\sum_{e\in\{-1,1\}}eN(v,w,l,e)
=
t_{w,l}
\sum_{e\in\{-1,1\}}eN(v,w,l+1,e).
\end{equation}
For each $e\in\{-1,1\}$, adding
and subtracting
\eqref{eq:alignment-count-total}
and
\eqref{eq:alignment-count-signed}
yields
\begin{equation}\label{eq:alignment-colored-matching}
N(v,w,l,e)
=
N(v,w,l+1,t_{w,l}e).
\end{equation}
Thus, for each auxiliary tower $v$,
each source tower $w$, and each
$l<h_w$, we can choose a bijection
between the auxiliary levels of
label $(w,l)$ and color $e$
and those of label $(w,l+1)$
and color $t_{w,l}e$.

Following these bijections partitions
all levels of the auxiliary tower into
disjoint chains
\[
(w,1),(w,2),\ldots,(w,h_w).
\]
Within each chain the colors satisfy
\[
e_{l+1}=t_{w,l}e_l.
\]
Every chain begins at a source base
label and ends at a source roof label. No level is omitted now, since every
nonroof level is matched to exactly
one level at the next source height,
and every nonbase level has exactly
one predecessor. No cycle can occur,
because the source height strictly
increases along each matching.

\smallskip
\noindent
\textbf{Step 4. Ordering the chains.}
Fix an auxiliary tower $v$.
A chain belonging to source tower
$w$ has first and last colors
$e_{\mathrm{in}}$ and $e_{\mathrm{out}}$
satisfying
\[
e_{\mathrm{out}}
=
\left(\prod_{l=1}^{h_w-1}t_{w,l}\right)
e_{\mathrm{in}}
=
\frac{T_w}{\tau_0}e_{\mathrm{in}}.
\]
Hence
\[
\tau_0e_{\mathrm{in}}
=
T_we_{\mathrm{out}}.
\]
Associate to this chain a directed
edge on the two-state set
$\{-1,1\}$, from
$T_we_{\mathrm{out}}$ to
$e_{\mathrm{out}}$.
If $T_w=1$, this edge is a loop.
If $T_w=-1$, it joins the two
different states.
Let
\[
n_+(v)
=
\#\{\text{odd chains ending in color }1\},
\]
and
\[
n_-(v)
=
\#\{\text{odd chains ending in color }-1\}.
\]
Each odd chain contributes exactly
one source roof level to $g$.
Therefore the last equality in
\eqref{eq:alignment-exact-sums}
gives
\begin{equation}\label{eq:alignment-odd-balance}
n_+(v)-n_-(v)
=
\frac{1-S_v}{2}.
\end{equation}
By the choice of the auxiliary tower
heights, every auxiliary tower contains
a level labeled by the roof of the
fixed odd source tower $w_*$.
Thus
\[
n_+(v)+n_-(v)>0.
\]

We now order all odd edges into
a directed path.

If $S_v=1$, then
\eqref{eq:alignment-odd-balance}
gives $n_+(v)=n_-(v)$.
Both numbers are positive.
Alternating edges
\[
1\longrightarrow-1,
\qquad
-1\longrightarrow1
\]
produces a path beginning and
ending at $1$, using every odd edge.

If $S_v=-1$, then
\[
n_+(v)=n_-(v)+1.
\]
Alternating the odd edges, beginning
with an edge from $-1$ to $1$,
produces a path from $-1$ to $1$
using every odd edge.

In either case the path starts at
$S_v$, ends at $1$, and visits
both states.
Every even chain corresponds to
a loop at one of these states.
Insert each such loop at an
occurrence of its state on the path.
This produces an ordering of all
chains, with no repetition or omission.

Within each chain retain the order
of its source levels.
We thereby obtain an ordering
\[
p_v(1),p_v(2),\ldots,p_v(r_v)
\]
of the original ambient levels.
By construction, the resulting
sequence has the following properties:
\begin{equation}\label{eq:alignment-order-properties}
\begin{aligned}
e_{v,p_v(r_v)}&=1,\\
e_{v,p_v(1)}&=\tau_0S_v.
\end{aligned}
\end{equation}
Moreover, if two consecutive
positions belong to the same
source chain, their labels are
consecutive source levels and
their colors satisfy the
corresponding source sign relation.

Finally, if one chain ends and another
begins, the first position of
the new chain has color $\tau_0$
times the last color of the
preceding chain.

\smallskip
\noindent
\textbf{Step 5. Realization by an element in the
topological full-group.}
Define $\sigma_1$ on the 
auxiliary level $Z_{v,p_v(j)}$ by
\[
\sigma_1|_{Z_{v,p_v(j)}}
=
\alpha^{j-p_v(j)}.
\]
The map sends that level onto
$Z_{v,j}$.

Since the levels form a finite
clopen partition and the assignments
permute the levels within each
ambient tower, $\sigma_1$ is a
homeomorphism belonging to
$[[\alpha]]$. Let 
\[\Sigma_1\]
be its iterate
lift, and correct the coefficient
homeomorphism by $\Sigma_1^{-1}$.

The orientation sign of the iterate
from old position $p$ to new
position $j$ is
\[
b_{v,p}b_{v,j}.
\]
This follows directly by multiplying
the signs along the orbit segment;
the same formula holds for positive
and negative displacements because
the signs belong to $\{-1,1\}$.
Consequently, the corrected
coefficient degree at new position
$j$ is
\[
\begin{aligned}
a'|_{Z_{v,j}}
&=
a|_{Z_{v,p_v(j)}}
b_{v,p_v(j)}b_{v,j}\\
&=
e_{v,p_v(j)}b_{v,j}.
\end{aligned}
\]
Consider first an interior level
$Z_{v,j}$, where $j<r_v$.
Since
\[
b_{v,j}=s_{v,j}b_{v,j+1},
\]
we obtain
\[
\begin{aligned}
a'|_{Z_{v,j+1}}\,
s_{v,j}\,
a'|_{Z_{v,j}}
&=
e_{v,p_v(j+1)}
b_{v,j+1}
s_{v,j}
e_{v,p_v(j)}
b_{v,j}\\
&=
e_{v,p_v(j+1)}
e_{v,p_v(j)}.
\end{aligned}
\]

If positions $j$ and $j+1$
belong to the same source chain,
their colors differ by the sign
$t_{w,l}$ of the current source
label. Therefore
\[
e_{v,p_v(j+1)}
e_{v,p_v(j)}
=
t_{w,l}.
\]
If $j$ is the last position of
one chain, its source label is
a roof label. The next position
is the first level of another
chain, and the chain ordering
gives
\[
e_{v,p_v(j+1)}
=
\tau_0e_{v,p_v(j)}.
\]
Hence
\[
e_{v,p_v(j+1)}
e_{v,p_v(j)}
=
\tau_0,
\]
which is precisely the source
sign on the merged roof.

It remains to verify the sign
equation across ambient towers.
By
\eqref{eq:alignment-order-properties},
the last position of every
ambient tower has color $1$.
Since $b_{v,r_v}=1$, its corrected
coefficient degree is also $1$.
Suppose that the next point
belongs to the base of ambient
tower $v'$. At that base,
\[
e_{v',p_{v'}(1)}
=
\tau_0S_{v'},
\qquad
b_{v',1}
=
\frac{S_{v'}}c.
\]
Thus
\[
a'|_{Z_{v',1}}
=
\tau_0S_{v'}\frac{S_{v'}}c
=
\frac{\tau_0}c.
\]
Since the sign on the preceding
auxiliary roof is $c$, we have
\[
a'(\alpha\omega)s(\omega)a'(\omega)
=
\frac{\tau_0}c\cdot c\cdot1
=
\tau_0.
\]
This is again the source sign
on the merged roof.
We have therefore proved, at
every point $\omega\in\Omega_1$,
\begin{equation}\label{eq:alignment-final-sign}
a'(\alpha\omega)s(\omega)a'(\omega)
=
t(h'(\omega)),
\end{equation}
where $h'$ is the base map of
the final homeomorphism of coefficient algebras.

We now verify the base identities.
The ordering preserves every
source chain and places its
successive labels consecutively.
Thus, for each nonroof source
label $W_{w,l}$,
\[
h'\alpha(h')^{-1}(W_{w,l})
\subset W_{w,l+1}.
\]
At a source roof label, the next
position belongs to a source base
label, so
\[
h'\alpha(h')^{-1}(R)
\subset\bigcup_wW_{w,1}.
\]
These inclusions are equalities.
Indeed, the source nonroof labels
and the merged roof partition
$\Omega_2$, and their prescribed
successors also partition
$\Omega_2$. Since
$h'\alpha(h')^{-1}$ is a
homeomorphism, the images of
the partition members are
pairwise disjoint and cover
$\Omega_2$.
Consequently,
\[
h'\alpha(h')^{-1}(U)
=
\beta(U)
\qquad(U\in\mathcal Q^\flat).
\]
Equation~\eqref{eq:alignment-final-sign}
is exactly
\[
\deg(H'\taf(H')^{-1})_\eta
=
t(\eta)
\qquad(\eta\in\Omega_2).
\]

Finally, the complete correction
is the composition of $\sigma_1$
with the finitely many full-group
corrections of Step 1.
Their iterate lifts compose
to the iterate lift of the
corresponding full-group element.
By Lemma~\ref{lem:full-group-inner-lift},
there is a unitary $W\in A$ such that
\[
j_{\taf}\xi'
=
\Ad(W)j_{\taf}\xi.
\]
This proves all the assertions.
\end{proof}

\begin{rem}
Proposition~\ref{prop:finite-signed-alignment} provides a
finite-level normalization of both the base dynamics and the
fiber orientation. Starting from a coefficient isomorphism
satisfying the prescribed $K$-theoretic compatibility, it
produces a correction that agrees with the target base
dynamics on a prescribed finite partition and matches the
fiber orientation signs exactly.

The essential point is that these two requirements can be
achieved simultaneously. The signed coinvariant relations
provide the necessary compatibility, while the
Kakutani--Rokhlin tower construction realizes it by
rearranging finitely many clopen levels. Since the
correction is implemented by an inner automorphism of the
crossed product, the induced coefficient $K$-theory maps
remain unchanged. This will be the key ingredient in
constructing approximate conjugacies with the required
$K$-theoretic control.
\end{rem}

\begin{lem}\label{lem:positive-global-sign}
Suppose that $s=t=1$, and let
$k=(k_0,k_1)$ be a graded isomorphism for which
the finite coefficient $K$-condition holds in both
directions. 

Let $H_n:Y\to X$ and $G_n:X\to Y$
be the homeomorphisms inducing the corresponding
coefficient isomorphisms. Then there exist
$\epsilon\in\{-1,1\}$ and $N\in\mathbb N$ such that
\[
\deg((H_n)_\eta)=\epsilon
\quad(\eta\in\Omega_2),
\qquad
\deg((G_n)_\omega)=\epsilon
\quad(\omega\in\Omega_1)
\]
for every $n\geq N$.
\end{lem}

\begin{proof}
Recall that
\[
D_\alpha
=
C(\Omega_1,\mathbb Z)/
(I-U_\alpha)C(\Omega_1,\mathbb Z),
\qquad
D_\beta
=
C(\Omega_2,\mathbb Z)/
(I-U_\beta)C(\Omega_2,\mathbb Z).
\]
We equip these groups with their natural positive
cones, consisting of classes represented by
nonnegative integer-valued continuous functions.

Recall that both groups are torsion-free by
Proposition~\ref{prop:signed-torsion}. Moreover,
every nonzero positive element is an order
unit. Indeed, every invariant probability measure
of a minimal Cantor system has full support.
Consequently, a nonnegative function representing
the zero coinvariant class must vanish identically.
On the other hand, if $f\geq 0$ is nonzero, then
minimality gives an integer $m\geq 1$ such that
\[
\sum_{j=0}^{m-1}f\circ\alpha^j\geq 1.
\]
Hence $m[f]\geq[1]$, so $[f]$ is an order unit.
The same argument applies to $D_\beta$.

Since $s=t=1$, the coefficient images in both
$K$-groups are identified with the corresponding
ordinary coinvariant groups. The finite coefficient
$K$-condition in both directions implies that
the restrictions of $k_0$ and $k_1$ induce group
isomorphisms
\[
K,L:D_\alpha\longrightarrow D_\beta,
\]
whose inverses are induced by $k_0^{-1}$ and
$k_1^{-1}$, respectively.

Write
\[
H_n(\eta,z)
=
(h_n(\eta),(H_n)_\eta(z)),
\]
and set
\[
a_n(\eta)=\deg((H_n)_\eta)
\in\{-1,1\}.
\]
The coefficient isomorphism induced by $H_n$
acts on fiber rank and fiber winding by
\[
f\longmapsto f\circ h_n,
\qquad
f\longmapsto
(a_n)(f\circ h_n),
\]
respectively.
Thus, for every $f\in C(\Omega_1,\mathbb Z)$,
the finite coefficient $K$-condition gives,
for all sufficiently large $n$,
\[
K[f]=[f\circ h_n],
\]
and
\[
L[f]
=
[(a_n)(f\circ h_n)].
\]
In particular,
\[
\begin{aligned}
(K+L)[f]
&=
2[1_{\{a_n=1\}}
(f\circ h_n)],\\
(K-L)[f]
&=
2[1_{\{a_n=-1\}}
(f\circ h_n)].
\end{aligned}
\]
Since $D_\beta$ is torsion-free, these identities
define unique homomorphisms
\[
P_+=\frac{K+L}{2},
\qquad
P_-=\frac{K-L}{2}
\]
from $D_\alpha$ to $D_\beta$. More precisely, the identities show
that $(K+L)[f]$ and $(K-L)[f]$ belong to
$2D_\beta$ for every $[f]\in D_\alpha$.
Injectivity of multiplication by $2$ guarantees
that the resulting elements are independent of
the representative and of the sufficiently large
index used to compute them. The same identities
show that both $P_+$ and $P_-$ are positive.

Applying the same argument to the reverse
coefficient sequence gives positive homomorphisms
\[
Q_+=\frac{K^{-1}+L^{-1}}{2},
\qquad
Q_-=\frac{K^{-1}-L^{-1}}{2}
\]
from $D_\beta$ to $D_\alpha$.

We now show that either $P_+=0$ or $P_-=0$.
By the definitions,
\[
4(Q_+P_-+Q_-P_+)
=
2(K^{-1}K-L^{-1}L)=0.
\]
Torsion-freeness of $D_\alpha$ gives
\[
Q_+P_-+Q_-P_+=0.
\]
Both compositions are positive. Since the
positive cone of $D_\alpha$ is proper,
\[
Q_+P_-(x)=Q_-P_+(x)=0
\qquad(x\in D_\alpha^+).
\]
As this cone generates $D_\alpha$, we obtain
\[
Q_+P_-=Q_-P_+=0.
\]

We use the following elementary consequence
of the order-unit property: the composition
of two nonzero positive homomorphisms between
$D_\alpha$ and $D_\beta$, in either order,
cannot vanish.

Indeed, let $P:D_\alpha\to D_\beta$ and
$Q:D_\beta\to D_\alpha$ be nonzero positive
homomorphisms. Since $[1]$ is an order unit
and the positive cones are generating,
\[
P([1])\neq 0,
\qquad
Q([1])\neq 0.
\]
The nonzero positive element $P([1])$ is
an order unit in $D_\beta$. Hence there is
$m\geq1$ such that
\[
mP([1])\geq[1].
\]
Applying $Q$ yields
\[
mQP([1])\geq Q([1])>0.
\]
Thus $QP\neq0$.

Suppose now that both $P_+$ and $P_-$ are
nonzero. The preceding observation, together
with
\[
Q_+P_-=Q_-P_+=0,
\]
forces $Q_+=Q_-=0$. This contradicts
\[
Q_++Q_-=K^{-1}.
\]
Therefore at least one of $P_+$ and $P_-$
vanishes. Since
\[
P_++P_-=K,
\qquad
P_+-P_-=L,
\]
there exists $\epsilon\in\{-1,1\}$ such that
\[
L=\epsilon K.
\]

We next show that this equality forces the
fiber degree of $H_n$ to be identically
$\epsilon$ for all sufficiently large $n$.
Since the coefficient isomorphisms are unital,
\[
K[1]=[1].
\]
Applying the finite coefficient $K$-condition
to the constant winding function $1$ gives
\[
[a_n]
=
L[1]
=
\epsilon[1]
\]
for all sufficiently large $n$.
If $\epsilon=1$, this equality implies
\[
2[1_{\{a_n=-1\}}]=0;
\]
if $\epsilon=-1$, it implies
\[
2[1_{\{a_n=1\}}]=0.
\]
In either case, torsion-freeness gives
\[
[1_{\{a_n=-\epsilon\}}]=0.
\]
Since a nonnegative function representing
the zero class must vanish,
\[
a_n\equiv\epsilon,
\]
and hence
\[
\deg((H_n)_\eta)=\epsilon
\qquad(\eta\in\Omega_2)
\]
for all sufficiently large $n$.

Finally, the reverse coefficient sequence
induces the isomorphisms
\[
K^{-1}:D_\beta\to D_\alpha,
\qquad
L^{-1}=\epsilon K^{-1}.
\]
Repeating the preceding argument with the
roles of $\alpha$ and $\beta$ interchanged
shows that
\[
\deg((G_n)_\omega)=\epsilon
\qquad(\omega\in\Omega_1)
\]
for all sufficiently large $n$. Taking a common $N$ for the two
sequences completes the proof.
\end{proof}

\begin{prop}\label{prop:simultaneous-normalization}
Under condition~\textup{(2)} of Theorem~\ref{thm:main}, there are
homeomorphisms $H_n:Y\to X$ and $K_n:X\to Y$, with base maps
$h_n:\Omega_2\to\Omega_1$ and $\ell_n:\Omega_1\to\Omega_2$, such that
\begin{equation}\label{eq:normalized-bases}
h_n\beta h_n^{-1}\longrightarrow\alpha,\qquad
\ell_n\alpha\ell_n^{-1}\longrightarrow\beta
\end{equation}
uniformly, and, eventually at every base point,
\begin{equation}\label{eq:normalized-signs}
\deg(H_n\tbt H_n^{-1})_\omega=s(\omega),\qquad
\deg(K_n\taf K_n^{-1})_\eta=t(\eta).
\end{equation}
At every index the pullbacks retain the entire embedded maps:
\begin{equation}\label{eq:normalized-coefficient-K}
(j_{\tbt}H_n^*)_{*i}=(j_{\tbt}\chi_n)_{*i},\qquad
(j_{\taf}K_n^*)_{*i}=(j_{\taf}\eta_n)_{*i},\quad i=0,1.
\end{equation}
\end{prop}

\begin{proof}
By Corollary~\ref{cor:orientation-K1-invariant},
the orientation signs $s$ and $t$ are either both trivial or
both nontrivial. We treat the two
possibilities separately.

\smallskip
\noindent
{\bf Case 1. Neither $s$ nor $t$ is a sign
coboundary.} Choose a sequence of first-return
Kakutani--Rokhlin partitions
\[
\mathcal Q_m
=
\{W_{w,l}^{(m)}:
1\leq w\leq q_m,\ 
1\leq l\leq h_w^{(m)}\}
\]
for $\beta$ such that $t$ is constant on
every level, the roof union $R_m$ is
contained in a clopen set on which $t$
is constant, and
\[
\operatorname{mesh}(\mathcal Q_m^\flat)
\longrightarrow 0,
\qquad
\operatorname{mesh}(\beta(\mathcal Q_m^\flat))
\longrightarrow 0,
\]
where $\mathcal Q_m^\flat$ consists of
the nonroof levels together with the
merged roof $R_m$. Such a sequence exists. Indeed, fix
$\eta_0\in\Omega_2$ and choose a decreasing
sequence of clopen neighborhoods $R_m$ of
$\eta_0$ whose diameters tend to zero and
on which $t$ is constant. Take the
first-return partitions over $\beta(R_m)$,
so that their merged roof unions are
precisely $R_m$. By refining the tower bases according to
finite clopen itineraries, we may arrange
that $t$ is constant on every level and
that the diameters of all nonroof levels
tend to zero. Consequently,
$\operatorname{mesh}(\mathcal Q_m^\flat)
\to 0$. Since $\beta$ is uniformly
continuous, we also have
$\operatorname{mesh}
(\beta(\mathcal Q_m^\flat))\to 0$.

For each $m$, let $E_m$ be the subgroup
of $C(\Omega_2,\mathbb Z)$ generated by
the characteristic functions of the
levels of $\mathcal Q_m$, which is a finitely
generated subgroup. By the finite coefficient $K$-condition,
there exists $N_m$ such that, for every
$n\geq N_m$ and $i=0,1$,
\[
(k_i^{-1}(j_{\tbt})_{*i})(f)
=
(j_{\taf}\eta_n)_{*i}(f)
\qquad(f\in E_m).
\]
Choose an increasing sequence
$m(n)\to\infty$ sufficiently slowly so that
$n\geq N_{m(n)}$ for all sufficiently
large $n$. For each such $n$, apply
Proposition~\ref{prop:finite-signed-alignment}
to $\eta_n$, with
\[
\lambda=k^{-1},
\qquad
\mathcal Q=\mathcal Q_{m(n)}.
\]
It gives a corrected coefficient
homeomorphism
\[
K_n:X\longrightarrow Y
\]
with base map
$\ell_n:\Omega_1\to\Omega_2$ satisfying
\[
\ell_n\alpha\ell_n^{-1}(U)
=
\beta(U)
\qquad
(U\in\mathcal Q_{m(n)}^\flat),
\]
and
\[
\deg(K_n\taf K_n^{-1})_\eta
=
t(\eta)
\qquad(\eta\in\Omega_2).
\]
Put
\[
\beta_n=\ell_n\alpha\ell_n^{-1}.
\]
For every $U\in\mathcal Q_{m(n)}^\flat$,
the preceding equality gives
$\beta_n(U)=\beta(U)$. Hence, for
$\eta\in U$,
\[
d(\beta_n(\eta),\beta(\eta))
\leq
\operatorname{diam}(\beta(U)).
\]
It follows that
\[
\sup_{\eta\in\Omega_2}
d(\beta_n(\eta),\beta(\eta))
\leq
\operatorname{mesh}
\bigl(\beta(\mathcal Q_{m(n)}^\flat)\bigr)
\longrightarrow0.
\]
Likewise, for $\eta\in\beta(U)$,
both $\beta_n^{-1}(\eta)$ and
$\beta^{-1}(\eta)$ belong to $U$.
Consequently,
\[
\sup_{\eta\in\Omega_2}
d(\beta_n^{-1}(\eta),\beta^{-1}(\eta))
\leq
\operatorname{mesh}
(\mathcal Q_{m(n)}^\flat)
\longrightarrow0.
\]
Thus
\[
\ell_n\alpha\ell_n^{-1}
\longrightarrow\beta
\]
uniformly, together with the inverses.
Moreover,
Proposition~\ref{prop:finite-signed-alignment}
and Lemma~\ref{lem:full-group-inner-lift}
give a unitary $V_n\in A$ such that
\[
j_{\taf}K_n^*
=
\operatorname{Ad}(V_n)
j_{\taf}\eta_n.
\]
Since inner automorphisms induce the
identity on $K$-theory,
\[
(j_{\taf}K_n^*)_{*i}
=
(j_{\taf}\eta_n)_{*i},
\qquad i=0,1.
\]

For the  other direction, interchange
the roles of the two systems and apply
the same construction to $\chi_n$,
using the graded isomorphism $k$.
This gives homeomorphisms
\[
H_n:Y\longrightarrow X
\]
with base maps
$h_n:\Omega_2\to\Omega_1$ such that
\[
h_n\beta h_n^{-1}\longrightarrow\alpha,
\]
\[
\deg(H_n\tbt H_n^{-1})_\omega
=
s(\omega)
\qquad(\omega\in\Omega_1),
\]
and
\[
(j_{\tbt}H_n^*)_{*i}
=
(j_{\tbt}\chi_n)_{*i},
\qquad i=0,1.
\]
Taking a common $N$ establishes
all the required assertions in this
case.

\smallskip
\noindent
{\bf Case 2. Both $s$ and $t$ are sign
coboundaries.} Choose continuous functions
\[
a:\Omega_1\to\{-1,1\},
\qquad
b:\Omega_2\to\{-1,1\}
\]
such that
\[
s(\omega)=a(\alpha\omega)a(\omega),
\qquad
t(\eta)=b(\beta\eta)b(\eta).
\]
Define the base-preserving reflection
homeomorphisms
\[
C_a(\omega,z)=(\omega,a(\omega)z),
\qquad
C_b(\eta,z)=(\eta,b(\eta)z),
\]
and put
\[
\taf^+=C_a\taf C_a^{-1},
\qquad
\tbt^+=C_b\tbt C_b^{-1}.
\]
By Lemma~\ref{lem:reflection-coordinates},
both $\taf^+$ and $\tbt^+$ have
orientation-preserving fiber maps.
Let
\[
\theta_a(f)=f\circ C_a^{-1},
\qquad
\theta_b(g)=g\circ C_b^{-1}
\]
be the corresponding coefficient
isomorphisms, and let
\[
\Theta_a:A\longrightarrow A^+,
\qquad
\Theta_b:B\longrightarrow B^+
\]
be their crossed-product extensions.
They satisfy the exact coefficient
identities
\[
\Theta_a j_{\taf}
=
j_{\taf^+}\theta_a,
\qquad
\Theta_b j_{\tbt}
=
j_{\tbt^+}\theta_b.
\]
Rewrite the prescribed graded
isomorphism and coefficient sequences
by setting
\[
k^+
=
(\Theta_b)_*\,k\,(\Theta_a)_*^{-1},
\]
\[
\chi_n^+
=
\theta_b\chi_n\theta_a^{-1},
\qquad
\eta_n^+
=
\theta_a\eta_n\theta_b^{-1}.
\]
The coefficient identities show
that these new objects again satisfy
the finite coefficient $K$-condition
in both directions.

Indeed, the inverse image of any
finitely generated coefficient
subgroup under either fixed
isomorphism of coefficient algebras is again
finitely generated. Applying the
original compatibility condition
to these inverse images
gives the required equalities in
both degrees.

Lemma~\ref{lem:positive-global-sign}
therefore provides
$\epsilon\in\{-1,1\}$ and $N_0$
such that, for every $n\geq N_0$,
the coefficient homeomorphisms
associated with $\chi_n^+$ and
$\eta_n^+$ have constant fiber
degree $\epsilon$ at every base
point.

Apply
Proposition~\ref{prop:combined-normalization}
to the new coefficient
sequences. We obtain corrected
homeomorphisms
\[
H_n^+:Y\longrightarrow X,
\qquad
K_n^+:X\longrightarrow Y
\]
whose base maps $h_n$ and $\ell_n$
satisfy
\[
h_n\beta h_n^{-1}\longrightarrow\alpha,
\qquad
\ell_n\alpha\ell_n^{-1}\longrightarrow\beta.
\]

The corrections supplied by
Proposition~\ref{prop:combined-normalization}
are iterate lifts of elements
of the corresponding topological
full groups. Since all fiber maps
of $\taf^+$ and $\tbt^+$ have degree
$1$, every such iterate lift has
fiber degree $1$. Consequently,
these corrections preserve the
constant fiber degree $\epsilon$.
Thus, for all sufficiently large $n$,
\[
\deg((H_n^+)_\eta)=\epsilon
\qquad(\eta\in\Omega_2),
\]
and
\[
\deg((K_n^+)_\omega)=\epsilon
\qquad(\omega\in\Omega_1).
\]
It follows that the fiber degrees
of the conjugated positive systems
are identically $1$:
\[
\deg(H_n^+\tbt^+(H_n^+)^{-1})_\omega
=
\epsilon^2=1,
\]
and
\[
\deg(K_n^+\taf^+(K_n^+)^{-1})_\eta
=
\epsilon^2=1.
\]

Now return to the original coordinates
by defining
\[
H_n=C_a^{-1}H_n^+C_b,
\qquad
K_n=C_b^{-1}K_n^+C_a.
\]
Since $C_a$ and $C_b$ fix the base
coordinates, the base maps of
$H_n$ and $K_n$ remain $h_n$
and $\ell_n$. Hence
\[
h_n\beta h_n^{-1}\longrightarrow\alpha,
\qquad
\ell_n\alpha\ell_n^{-1}\longrightarrow\beta.
\]
Put
\[
\alpha_n=h_n\beta h_n^{-1}.
\]
The coordinate identities give
\[
H_n\tbt H_n^{-1}
=
C_a^{-1}
\bigl(H_n^+\tbt^+(H_n^+)^{-1}\bigr)
C_a.
\]
Since the middle homeomorphism
has fiber degree $1$, we obtain
\[
\deg(H_n\tbt H_n^{-1})_\omega
=
a(\alpha_n\omega)a(\omega).
\]
Since $a$ is continuous and
takes values in the discrete set
$\{-1,1\}$, uniform convergence
$\alpha_n\to\alpha$ implies
\[
a(\alpha_n\omega)=a(\alpha\omega)
\]
for every $\omega\in\Omega_1$
and all sufficiently large $n$.
Therefore
\[
\deg(H_n\tbt H_n^{-1})_\omega
=
a(\alpha\omega)a(\omega)
=
s(\omega)
\]
eventually, uniformly over
$\omega\in\Omega_1$. Similarly, putting
\[
\beta_n=\ell_n\alpha\ell_n^{-1},
\]
we obtain
\[
\deg(K_n\taf K_n^{-1})_\eta
=
b(\beta_n\eta)b(\eta).
\]
Since $\beta_n\to\beta$ uniformly,
local constancy of $b$ gives
\[
b(\beta_n\eta)=b(\beta\eta)
\]
for every $\eta\in\Omega_2$
and all sufficiently large $n$.
Consequently,
\[
\deg(K_n\taf K_n^{-1})_\eta
=
t(\eta).
\]

It remains to verify preservation
of the
$K$-theory maps induced from canonical embeddings.
By
Proposition~\ref{prop:combined-normalization},
there are unitaries
$W_n^+\in B^+$ and $V_n^+\in A^+$
such that
\[
j_{\tbt^+}(H_n^+)^*
=
\operatorname{Ad}(W_n^+)
j_{\tbt^+}\chi_n^+,
\]
and
\[
j_{\taf^+}(K_n^+)^*
=
\operatorname{Ad}(V_n^+)
j_{\taf^+}\eta_n^+.
\]
Using the exact coefficient
identities for $\Theta_a$ and
$\Theta_b$, together with the
definitions of the new
coefficient maps, gives
\[
j_{\tbt}H_n^*
=
\operatorname{Ad}
\bigl(\Theta_b^{-1}(W_n^+)\bigr)
j_{\tbt}\chi_n
\]
and
\[
j_{\taf}K_n^*
=
\operatorname{Ad}
\bigl(\Theta_a^{-1}(V_n^+)\bigr)
j_{\taf}\eta_n.
\]
Hence, for $i=0,1$,
\[
(j_{\tbt}H_n^*)_{*i}
=
(j_{\tbt}\chi_n)_{*i},
\qquad
(j_{\taf}K_n^*)_{*i}
=
(j_{\taf}\eta_n)_{*i}.
\]
The preceding constructions apply
for all sufficiently large indices.
For the finitely many remaining
indices, retain the original
coefficient homeomorphisms. The
$K$-maps induced from the canonical embeddings are
then preserved at every index,
while the asymptotic base and
orientation conclusions remain
unchanged. This completes the proof.
\end{proof}

\begin{rem}\label{rem:normalization-scope}
We note that the normalization so far gives no
approximation of the fiber homeomorphisms themselves, which is exactly what we are going to consider in the next section.
\end{rem}

\section{Return-tower corrections}\label{sec:towers}
Proposition~\ref{prop:simultaneous-normalization}
provides homeomorphisms whose conjugated extensions
have base maps converging uniformly to the target
base maps and whose fiber degrees agree with the
target orientation signs. It remains to correct the
fiber homeomorphisms themselves.

The main difficulty is that a correction along a
return tower must match the entire return
homeomorphism while keeping every intermediate
fiber map close to the target. We overcome this
difficulty by combining endpoint interpolation with
a uniform width estimate derived from minimality.
The resulting corrections are isotopic to the
identity and hence preserve the entire coefficient
$K$-theory.

The construction requires only minimality of the
extensions. Rigidity will enter later, when we
recover the required trace estimates.

\subsection{Lifts and endpoint interpolation}
We begin with two useful Lemmas.
\begin{lem}\label{lem:continuous-lifts}
Let $\Omega$ be a Cantor space and let
$\omega\mapsto\varphi_\omega$ be a continuous
family of circle homeomorphisms. Then there
exists a jointly continuous family of real
lifts $p_\omega:\mathbb R\to\mathbb R$
satisfying
\[
p_\omega(u+1)=p_\omega(u)+s(\omega),
\qquad
s(\omega)=\deg(\varphi_\omega).
\]
The inverse family is also jointly continuous.

Moreover, every base-preserving homeomorphism
$H$ of $\Omega\times\T$ whose fiber maps
preserve orientation is isotopic to the identity
through base-preserving homeomorphisms.
In particular,
\[
(H^*)_{*i}=((H^{-1})^*)_{*i}
=\operatorname{id}_{K_i(C(\Omega\times\T))},
\qquad i=0,1.
\]
\end{lem}

\begin{proof}
Since the degree is locally constant, the sets
\[
\Omega_\pm
=
\{\omega\in\Omega:s(\omega)=\pm1\}
\]
are clopen. We may therefore construct the
lifts separately on these two sets.

The continuous map
$\omega\mapsto\varphi_\omega(0)$
admits local continuous real lifts. By
compactness and zero-dimensionality of
$\Omega$, a finite clopen refinement gives
a continuous function $r:\Omega\to\mathbb R$
such that
\[
r(\omega)\pmod{\mathbb Z}
=
\varphi_\omega(0).
\]
For each $\omega$, let $p_\omega$ be the
unique lift of $\varphi_\omega$ satisfying
$p_\omega(0)=r(\omega)$. Its monotonicity
is determined by $s(\omega)$, and
\[
p_\omega(u+1)=p_\omega(u)+s(\omega).
\]

To see joint continuity, fix $\omega_0$.
For $\omega$ sufficiently close to
$\omega_0$, the circle-valued difference
$\varphi_\omega-\varphi_{\omega_0}$
has a unique continuous real representative
with values in $(-1/2,1/2)$.
Adding this representative to
$p_{\omega_0}$ gives a continuous local
family of lifts. The normalization at zero
identifies these lifts with $p_\omega$
after shrinking the neighborhood if
necessary. Thus $(\omega,u)\mapsto
p_\omega(u)$ is jointly continuous.
The same holds for the inverse family
by strict monotonicity and continuity.

Now let $H$ be a base-preserving
homeomorphism whose fiber maps preserve
orientation, and choose jointly continuous
increasing lifts $q_\omega$.
For $v\in[0,1]$, put
\[
q_{\omega,v}(u)
=
(1-v)u+vq_\omega(u).
\]
Each $q_{\omega,v}$ is a strictly
increasing homeomorphism of $\mathbb R$
satisfying
\[
q_{\omega,v}(u+1)=q_{\omega,v}(u)+1.
\]
Consequently,
\[
H_v(\omega,\bar u)
=
\bigl(\omega,\overline{q_{\omega,v}(u)}\bigr)
\]
defines an isotopy through
base-preserving homeomorphisms from
the identity to $H$.
Homotopy invariance of $K$-theory gives
\[
(H^*)_{*i}
=
((H^{-1})^*)_{*i}
=
\operatorname{id}_{K_i(C(\Omega\times\T))},
\qquad i=0,1.
\]
This completes the proof.
\end{proof}

\begin{lem}\label{lem:inverse-convergence}
If homeomorphisms $S_n$ of a compact metric space converge uniformly
to a homeomorphism $S$, then $S_n^{-1}\to S^{-1}$ uniformly.
\end{lem}
\begin{proof}
Uniformly in $z$,
$d(S(S_n^{-1}z),z)\le\sup_w d(Sw,S_nw)$.
Apply uniform continuity of $S^{-1}$.
\end{proof}

For any two increasing degree-one homeomorphisms
$u,v:\mathbb R\to\mathbb R$, the periodic function
$u-v$ has oscillation at most $1$. Indeed, after
translating one argument by an integer, we may
assume $t\le r\le t+1$, so that both
$u(r)-u(t)$ and $v(r)-v(t)$ belong to $[0,1]$.

Consequently, suppose that $A^-,A^+$ and $D_0$
are continuous families of increasing degree-one
homeomorphisms over a compact clopen subspace
$E$ of a Cantor space, and that
\[
A^+(x,u)-A^-(x,u)>4
\]
for all $(x,u)\in E\times\mathbb R$.
A finite clopen refinement allows us to choose
a locally constant function $k:E\to\mathbb Z$
such that
\[
1<D_0(x,0)+k(x)-A^-(x,0)<3.
\]
The oscillation estimate then shows that
$D(x,u)=D_0(x,u)+k(x)$ satisfies
\begin{equation}\label{eq:return-lift-bracket}
A^-(x,u)<D(x,u)<A^+(x,u)
\end{equation}
for all $(x,u)\in E\times\mathbb R$.

\begin{prop}
\label{prop:endpoint-interpolation}
Let $E$ be a compact metrizable space, let
$m\geq1$, and let $\delta>0$. Suppose that
$D$ and $a_j$, $0\leq j<m$, are continuous
families over $E$ of increasing degree-one
homeomorphisms of $\mathbb R$.

For $\lambda\in[-\delta,\delta]$, define
\[
A_0(x,\lambda,u)=u,
\qquad
A_{j+1}(x,\lambda,u)
=
a_j(x,A_j(x,\lambda,u))+\lambda.
\]
Assume that
\[
A_m(x,-\delta,u)
<
D(x,u)
<
A_m(x,\delta,u)
\]
for every $(x,u)\in E\times\mathbb R$.

Then there are continuous families
$b_j$, $0\leq j<m$, of increasing
degree-one homeomorphisms such that
\[
\sup_{x,u}|b_j(x,u)-a_j(x,u)|<\delta
\]
and
\[
b_{m-1}(x,\cdot)\circ\cdots\circ
b_0(x,\cdot)=D(x,\cdot).
\]
The inverse families are also jointly
continuous.
\end{prop}

\begin{proof}
For each $(x,u)$, the map
$\lambda\mapsto A_m(x,\lambda,u)$
is continuous and strictly increasing.
The strict endpoint inequalities therefore
determine a unique
$c(x,u)\in(-\delta,\delta)$ satisfying
\begin{equation}\label{eq:endpoint-parameter}
A_m(x,c(x,u),u)=D(x,u).
\end{equation}
Joint continuity and uniqueness imply
that $c$ is continuous. Since all the
maps involved commute with integer
translations, uniqueness also gives
\[
c(x,u+1)=c(x,u).
\]
Compactness of $E\times[0,1]$ yields
$\delta_0<\delta$ such that
\[
|c(x,u)|\leq\delta_0
\]
for every $(x,u)$.

Put
\[
B_j(x,u)=A_j(x,c(x,u),u).
\]
Then $B_0(x,u)=u$, $B_m=D$, and
\begin{equation}\label{eq:endpoint-recurrence}
B_{j+1}(x,u)
=
a_j(x,B_j(x,u))+c(x,u).
\end{equation}
Each $B_j$ is continuous and satisfies
\[
B_j(x,u+1)=B_j(x,u)+1.
\]

We claim that $B_j(x,\cdot)$ is strictly
increasing for every $j$. Fix $x$ and
$u<v$. If $c(x,u)\leq c(x,v)$, the
recurrence immediately gives
\[
B_j(x,u)<B_j(x,v)
\]
for every $j$. Otherwise, if this
inequality failed at some $j$, the strict inequality
$c(x,u)>c(x,v)$ would propagate the
reversed order through all subsequent
indices. This would imply
\[
D(x,u)>D(x,v),
\]
contrary to the strict monotonicity of
$D(x,\cdot)$. Thus every $B_j(x,\cdot)$
is strictly increasing. Its degree-one
periodicity makes it a homeomorphism
of $\mathbb R$.

The inverse families $B_j^{-1}$ are
jointly continuous. In fact,
by joint continuity, compactness, and
degree-one periodicity,
\[
M_j:=\sup_{(x,u)\in E\times\mathbb R}
|B_j(x,u)-u|<\infty.
\]
Hence
\[
|B_j(x,\cdot)^{-1}(v)-v|\leq M_j.
\]
If $(x_n,v_n)\to(x,v)$ and
$u_n=B_j(x_n,\cdot)^{-1}(v_n)$,
then $(u_n)$ is bounded, and every
convergent subsequence satisfies
\[
B_j(x,\lim_{k\to\infty}u_{n_k})
=\lim_{k\to\infty}v_{n_k}=v.
\]
By injectivity, all such subsequences
converge to $B_j(x,\cdot)^{-1}(v)$.
Therefore the inverse family is jointly
continuous.

Now we define
\[
b_j(x,\cdot)
=
B_{j+1}(x,\cdot)\circ B_j(x,\cdot)^{-1}.
\]
These are continuous families of
increasing degree-one homeomorphisms,
with jointly continuous inverses.
By \eqref{eq:endpoint-recurrence},
\[
b_j(x,v)
=
a_j(x,v)
+
c(x,B_j(x,\cdot)^{-1}(v)).
\]
Consequently,
\[
\sup_{x,v}|b_j(x,v)-a_j(x,v)|
\leq\delta_0<\delta.
\]
Finally, the composition telescopes:
\[
b_{m-1}(x,\cdot)\circ\cdots\circ
b_0(x,\cdot)
=
B_m(x,\cdot)\circ B_0(x,\cdot)^{-1}
=
D(x,\cdot).
\]
The proof is complete.
\end{proof}

\subsection{Endpoint width from minimality}
The purpose of this subsection is to show that minimality provides sufficient freedom to correct the fiber maps along long return towers. More precisely, we shall prove that the separation between the endpoint maps obtained by adding small positive and negative perturbations grows uniformly with the length of the iteration. This separation, which we call the {\it endpoint width}, will allow us to apply Proposition \ref{prop:endpoint-interpolation}.

To handle orientation-reversing fiber maps, we first pass to normalized lifts on the orientation cover. Let $p_\omega$ be jointly continuous real lifts
of $\varphi_\omega$ with degree $s(\omega)$. On \(\widehat\Omega=\Omega\times\{-1,1\}\), define as before
\[
\widehat\alpha(\omega,e)
=(\alpha\omega,s(\omega)e),
\qquad
a_{\omega,e}(u)
=s(\omega)e\,p_\omega(eu).
\]
Each \(a_{\omega,e}\) is an increasing degree-one homeomorphism of \(\mathbb R\), regardless of the original orientation sign. 
For $\delta>0$, define
\[
a_{\omega,e}^{\pm\delta}(u)
=a_{\omega,e}(u)\pm\delta.
\]
Let
\[
A_{\pm\delta}^{(n)}(\omega,e,\cdot)
=
a_{\widehat\alpha^{n-1}(\omega,e)}^{\pm\delta}
\circ\cdots\circ
a_{\widehat\alpha(\omega,e)}^{\pm\delta}
\circ a_{\omega,e}^{\pm\delta},
\]
and set
\[
W_n(\omega,e,u)
=
A_\delta^{(n)}(\omega,e,u)
-
A_{-\delta}^{(n)}(\omega,e,u).
\]

\begin{lem}\label{lem:finite-crossing}
Assume that $\alpha$ is minimal and fix
$\delta>0$. Then
either
\[
W_n(\omega,e,u)\leq 2
\]
for every $n,\omega,e,u$, or there exists
$B\geq1$ such that
\[
W_n(\omega,e,u)\geq
\left\lfloor\frac{n}{B}\right\rfloor
\]
for all $n,\omega,e,u$.
The second alternative holds whenever
$W_N(\omega_0,e_0,u_0)>2$ for some
$N,\omega_0,e_0,u_0$.

The same dichotomy holds for the products
of $p\pm\delta$ associated with any continuous
family of increasing degree-one homeomorphisms
over a minimal homeomorphism of a compact
metrizable space.
\end{lem}

\begin{proof}
Suppose first that
$W_N(\omega_0,e_0,u_0)>2$.
The oscillation estimate for increasing
degree-one maps implies
\[
W_N(\omega_0,e_0,u)>1
\qquad (u\in\R).
\]
By continuity, there is a nonempty open
neighborhood $V$ of $\omega_0$ such that
\[
W_N(\omega,e_0,u)>1
\qquad (\omega\in V,\ u\in\R).
\]
The reflection identity
\[
a_{\omega,-e}(u)
=-a_{\omega,e}(-u)
\]
gives
\[
W_N(\omega,-e,u)
=W_N(\omega,e,-u),
\]
so the same estimate holds for both signs.

By minimality, there exists $R\geq0$
such that every forward $\alpha$-orbit
meets $V$ within $R$ steps.

Suppose that $u_+\geq u_-$ and
$u_+-u_-\geq k$ for some integer $k\geq0$.
Since the perturbed maps are increasing
and commute with integer translations,
\[
A_\delta^{(N)}(\omega,e,u_+)
-
A_{-\delta}^{(N)}(\omega,e,u_-)
\geq
W_N(\omega,e,u_-)+k.
\]
Thus, whenever $\omega\in V$, the
difference is greater than $k+1$.

By minimality, every orbit enters $V$
within $R$ steps. The integer lower bound
is preserved under further iterations,
and an $N$-step composition beginning
over $V$ increases this bound by at least
one. Hence every block of
$B=R+N$ consecutive steps increases the
integer lower bound by at least one.
Consequently,
\[
W_n(\omega,e,u)\geq
\left\lfloor\frac{n}{B}\right\rfloor
\]
for all $n,\omega,e,u$.

If no finite crossing larger than $2$ occurs,
the first alternative holds.

For an increasing degree-one cocycle,
the same argument applies without the
sign coordinate or the reflection
identity. Minimality of the base again
provides the uniform return bound $R$,
and the integer separation argument
is unchanged.
\end{proof}

We shall also use the following standard
properties for
increasing degree-one homeomorphisms without further proofs.
\begin{lem}\label{lem:monotone-translation}
For a nondecreasing map $v:\R\to\R$ with $v(u+1)=v(u)+1$, continuity
not being required, there is $r\in\R$ such that
$|v^n(u)-u-nr|\le1$ for all $n$. 
In particular, $v^q(u)\ge u+c$ for all $u$
implies $qr\ge c$, and the reverse inequalities give $qr\le c$.
\end{lem}

\begin{prop}
\label{prop:positive-width}
Let $Z$ be a compact metrizable space, and let
\[
S(x,\bar u)
=
(\beta x,\overline{p_x(u)})
\]
be a minimal homeomorphism on $Z\times\T$, where
$\beta\in\Homeo(Z)$ and $p_x$ is a
continuous family of increasing degree-one
homeomorphisms of $\R$.

Then for every $\epsilon>0$, there exists
$L\geq1$ such that
\[
P_\epsilon^{(n)}(x,u)
-
P_{-\epsilon}^{(n)}(x,u)>4
\]
for all $n\geq L$, $x\in Z$, and $u\in\R$, where $P_{\pm\epsilon}^{(n)}$ denotes
the $n$-fold cocycle composition of
the maps $p_x\pm\epsilon$ over $\beta$.
\end{prop}

\begin{proof}
By Lemma~\ref{lem:finite-crossing},
it suffices to exclude uniformly bounded
endpoint widths. Suppose, to the contrary,
that for some $C>0$,
\[
0\leq
P_\epsilon^{(n)}(x,u)
-
P_{-\epsilon}^{(n)}(x,u)
\leq C
\]
for all $n\geq1$, $x\in Z$, and $u\in\R$.
Define
\[
S_\pm(x,u)
=
(\beta x,p_x(u)\pm\epsilon),
\]
and let $h_{n,x}$ be the fiber map of
$S_+^nS_-^{-n}$, with $h_{0,x}(u)=u$.
These maps are continuous, increasing,
and of degree one, and satisfy
\[
u\leq h_{n,x}(u)\leq u+C.
\]
Since $S_+S_-^{-1}$ is translation by
$2\epsilon$ on each fiber,
\[
h_{n+1,x}(u)>h_{n,x}(u).
\]
Consequently,
\[
h_x(u)=\sup_{n\geq0}h_{n,x}(u)
\]
is finite, nondecreasing, jointly lower
semicontinuous, and satisfies
\[
h_x(u+1)=h_x(u)+1.
\]
The identity
\[
S_+^{n+1}S_-^{-(n+1)}
=
S_+(S_+^nS_-^{-n})S_-^{-1}
\]
gives
\[
h_{n+1,\beta x}(p_x(u)-\epsilon)
=
p_x(h_{n,x}(u))+\epsilon.
\]
Passing to the limit and using
monotonicity, we obtain
\[
h_{\beta x}(p_x(u))
\geq p_x(h_x(u))+\epsilon.
\]
Since $h_{\beta x}$ is nondecreasing,
an induction on $q$ gives
\begin{equation}\label{eq:width-obstruction-gap}
h_{\beta x}^{q}(p_x(u))
\geq p_x(h_x^{q}(u))+\epsilon
\end{equation}
for every $q\geq1$.
Indeed, if the inequality holds for $q$,
then
\[
\begin{aligned}
h_{\beta x}^{q+1}(p_x(u))
&\geq h_{\beta x}(p_x(h_x^q(u)))\geq p_x(h_x^{q+1}(u))+\epsilon.
\end{aligned}
\]

We now claim that $h_x^q(u)$ is jointly
lower semicontinuous. Indeed,
\[
h_x^q(u)
=
\sup_{n\geq0}h_{n,x}^q(u).
\]
To verify this identity, note that
$h_{n,x}^q(u)\leq h_x^q(u)$.
Conversely, for $n\geq k$,
\[
h_{n,x}^q(u)
\geq h_{k,x}(h_{n,x}^{q-1}(u)).
\]
Induction on $q$, followed by
continuity of $h_{k,x}$ and
taking the supremum over $k$,
gives the reverse inequality.
The claim follows.

For $m\in\Z$ and $q\geq1$, put
\[
E_{m,q}
=
\{(x,\bar u)\in Z\times\T:
h_x^q(u)\leq u+m\}.
\]
These sets are closed. Moreover,
if $S(x,\bar u)\in E_{m,q}$, then
\eqref{eq:width-obstruction-gap} gives
\[
p_x(h_x^q(u))+\epsilon
\leq p_x(u)+m
=p_x(u+m).
\]
Since $p_x$ is strictly increasing,
\[
h_x^q(u)<u+m.
\]
Thus
\[
S^{-1}(E_{m,q})\subseteq E_{m,q}.
\]
Minimality implies that every nonempty
$E_{m,q}$ equals $Z\times\T$. By uniform continuity, choose $\eta>0$
such that
\[
|v-w|\leq\eta
\quad\Longrightarrow\quad
|p_x(v)-p_x(w)|<\epsilon
\]
for all $x\in Z$ and $v,w\in\R$.
Whenever $E_{m,q}=Z\times\T$,
\eqref{eq:width-obstruction-gap} gives
\[
p_x(h_x^q(u))+\epsilon
\leq p_x(u+m),
\]
and hence
\begin{equation}\label{eq:width-uniform-gap}
h_x^q(u)\leq u+m-\eta
\end{equation}
for all $x,u$. Fix $x_0\in Z$, and let $r$ be the
translation number of $h_{x_0}$.
Choose integers $q\geq1$ and $m$
such that
\[
0\leq m-qr<\eta.
\]
Such integers exist by rationality
or one-sided approximation modulo one.

We now claim that $E_{m,q}$ is nonempty.
Otherwise, lower semicontinuity and
compactness give $c>0$ such that
\[
h_{x_0}^q(u)\geq u+m+c
\qquad (u\in\R).
\]
Lemma~\ref{lem:monotone-translation}
would imply
\[
qr\geq m+c,
\]
contrary to $qr\leq m$. Therefore, $E_{m,q}=Z\times\T$.
By \eqref{eq:width-uniform-gap} and
Lemma~\ref{lem:monotone-translation},
\[
qr\leq m-\eta,
\]
contradicting $m-qr<\eta$.

This contradiction shows that the
endpoint widths cannot be uniformly
bounded. Lemma~\ref{lem:finite-crossing}
therefore gives uniform growth, and
the assertion follows.
\end{proof}

\begin{prop}\label{prop:signed-width}
Let
\[
\zeta(\omega,\bar u)
=(\alpha\omega,\varphi_\omega(\bar u))
\]
be a minimal Cantor--circle extension over $\af$, and let
$s(\omega)=\deg(\varphi_\omega)$.
Then for every $\delta>0$, we have
\[W_n(\omega,e,u)>4\]
for all $\omega\in \Omega$, $e\in\{-1,1\}$,
$u\in\R$, and all sufficiently large $n$.
\end{prop}

\begin{proof}
Consider the orientation-preserving extension
\[
\mathcal S(\omega,e,\bar u)
=
\bigl(
\alpha\omega,
s(\omega)e,
\overline{a_{\omega,e}(u)}
\bigr)
\]
over $\widehat\alpha$, where
\[
\widehat\alpha(\omega,e)
=
(\alpha\omega,s(\omega)e).
\]
By construction, its fiber maps have degree one.

Suppose first that $s$ is not a continuous
sign coboundary. By
Lemma~\ref{lem:orientation-cover},
the pullback extension
\[
\widehat \zeta(\omega,e,\bar u)
=
\bigl(
\alpha\omega,
s(\omega)e,
\varphi_\omega(\bar u)
\bigr)
\]
is minimal. The homeomorphism
\[
J(\omega,e,\bar u)
=
(\omega,e,e\bar u)
\]
conjugates $\widehat \zeta$ to $\mathcal S$.
Thus $\mathcal S$ is minimal, and
Proposition~\ref{prop:positive-width}
gives the required width estimate.

Now suppose that
\[
s(\omega)=b(\alpha\omega)b(\omega)
\]
for some continuous
$b:\Omega\to\{-1,1\}$.
For each $c\in\{-1,1\}$, let
\[
B_c=
\{(\omega,e):e=cb(\omega)\}.
\]
The sets $B_1$ and $B_{-1}$ form a
clopen $\widehat\alpha$-invariant
partition of $\Omega\times\{-1,1\}$.
Moreover, the restriction of
$\mathcal S$ to $B_c\times\T$
is conjugate to $\zeta$ through
\[
(\omega,e,\bar u)
\longmapsto
(\omega,e\bar u).
\]
Consequently, both restrictions
are minimal. 
Applying Proposition~\ref{prop:positive-width}
to both components, we obtain
$W_n(\omega,e,u)>4$ for all
$\omega$, $e$, and $u$, whenever
$n$ is sufficiently large.
\end{proof}

\subsection{Repair over nearby base maps}

We now use the endpoint width obtained above to
correct circle extensions over base maps close to
a fixed minimal homeomorphism. The correction is
made on clopen return towers. On each tower,
the initial orbit segment is left unchanged, while
the final finitely many fiber maps are adjusted to
match the prescribed return map. 

In this way, the
dynamical approximation is achieved without
changing the $K$-theory maps. The
modified levels also have uniformly small measure
for all invariant probabilities on the base space.

\begin{prop}
\label{prop:nearby-base-repair}
Let
\[
\taf(\omega,z)
=
(\alpha\omega,\varphi_\omega(z))
\]
be a minimal Cantor--circle extension, and put
$s(\omega)=\deg(\varphi_\omega)$.
Suppose that $\alpha_n$ are minimal
homeomorphisms converging uniformly to $\alpha$,
and that
\[
S_n(\omega,z)
=
(\alpha_n\omega,\theta_{n,\omega}(z))
\]
are continuous circle extensions satisfying
\[
\deg(\theta_{n,\omega})=s(\omega)
\]
for every $\omega$ and all sufficiently large $n$.
Put
\[
T_n(\omega,z)
=
(\alpha_n\omega,\varphi_\omega(z)).
\]

For every $\delta,\eta>0$, there are an integer
$L\geq1$, finitely many clopen sets $V_r$ as first return bases for $\alpha$, with heights
\[
r>L,\qquad \frac{L}{r}<\eta,
\]
and an integer $N$ such that, for every $n\geq N$,
the sets
\[
Q_{r,j}
=
\alpha^j(V_r)
=
\alpha_n^j(V_r),
\qquad 0\leq j<r,
\]
form the same first-return tower partition
for $\alpha$ and $\alpha_n$.

Moreover, for each $n\geq N$ there is a
base-preserving homeomorphism $C_n$,
isotopic to the identity, with the following
properties:

\begin{enumerate}
\item
The circle coordinates of
$C_nS_nC_n^{-1}$ and $T_n$
are uniformly within $\delta$.

\item
The map $C_n$ is the identity on
$V_r\times\T$ for every $r$, and
\[
C_nS_n^j(z)=T_n^j(z)
\]
for every $z\in V_r\times\T$ and
$0\leq j\leq r-L$.
Furthermore, the fiber maps of
$C_nS_nC_n^{-1}$ and $T_n$
agree outside
\[
E=
\bigcup_r\bigcup_{j=r-L}^{r-1}Q_{r,j}.
\]

\item
For every
$\mu\in M(\alpha)\cup M(\alpha_n)$,
\[
\mu(E)<\eta.
\]

\item
The pullbacks induced by $C_n$ and
$C_n^{-1}$ satisfy
\[
(C_n^*)_{*i}
=
((C_n^{-1})^*)_{*i}
=
\operatorname{id}_{K_i(C(\Omega\times\T))},
\qquad i=0,1.
\]
\end{enumerate}

Consequently, the homeomorphisms $C_n$
may be chosen so that
\[
C_nS_nC_n^{-1}\longrightarrow\taf
\]
uniformly, together with their inverses.

If $S_n=h_nGh_n^{-1}$ for homeomorphisms
$h_n$ and a circle extension $G$, then
$k_n=C_nh_n$ satisfies
\[
k_nGk_n^{-1}\longrightarrow\taf
\]
uniformly, with uniform inverse convergence.
Moreover, $k_n$ and $h_n$ induce the same
maps on $K_i$, as do $k_n^{-1}$ and
$h_n^{-1}$, for $i=0,1$.
\end{prop}

\begin{proof}
Fix $\delta,\eta>0$, and choose continuous
signed lifts $p_\omega$ of $\varphi_\omega$
by Lemma~\ref{lem:continuous-lifts}.

Applying Proposition~\ref{prop:signed-width}
with this $\delta$, we obtain an integer
$L\geq1$ and $\varepsilon>0$ such that
\[
W_L(\omega,e,u)\geq4+\varepsilon
\]
for all $\omega\in\Omega$,
$e\in\{-1,1\}$, and $u\in\R$.
For every fixed $j$, uniform convergence
$\alpha_n\to\alpha$ implies
$\alpha_n^j\to\alpha^j$ uniformly.
Since $s$ is locally constant,
\[
s(\alpha_n^j\omega)
=
s(\alpha^j\omega)
\]
for all $\omega$, all $0\leq j<L$,
and all sufficiently large $n$.
Let $W_L^{(n)}$ denote the endpoint width
defined using the fixed target lifts
$p_\omega$, with $\alpha$ replaced by
$\alpha_n$.
Since only finitely many iterates are
involved, the corresponding endpoint
products converge uniformly to those
defining $W_L$.
Consequently, for all sufficiently large $n$,
\[
W_L^{(n)}(\omega,e,u)>4
\]
for every $\omega\in\Omega$,
$e\in\{-1,1\}$, and $u\in\R$.

Choose a clopen set for
$\alpha$ whose return times satisfy
\[
r>L,\qquad L/r<\eta,
\]
and let $V_r$ be its nonempty
first-return bases. The sets
\[
Q_{r,j}=\alpha^j(V_r),
\qquad 0\leq j<r,
\]
form a finite clopen partition of
$\Omega$.
For each member $Q$ of this partition,
uniform convergence gives
\[
\alpha_n(Q)\subseteq\alpha(Q)
\]
for all sufficiently large $n$.
Since both families of images form
partitions of $\Omega$, these inclusions
are equalities. It follows that
\[
\alpha_n^j(V_r)=\alpha^j(V_r),
\qquad 0\leq j\leq r.
\]
In particular, for all sufficiently large $n$,
the sets $V_r$ have the same first-return times
under $\alpha_n$ as under $\alpha$, and
\[
\alpha_n^j(V_r)=\alpha^j(V_r),
\qquad 0\leq j\leq r.
\]
Increasing $N$ if necessary, we may assume
that the preceding tower identities and
width estimates, together with
\[
\deg(\theta_{n,\omega})=s(\omega),
\]
hold for every $n\geq N$ and
$\omega\in\Omega$.

Fix such an $n$ and one of the towers.
For $\omega\in V_r$, write
\[
\omega_j=\alpha_n^j\omega,
\qquad
e_0=1,\qquad
e_{j+1}=s(\omega_j)e_j,
\]
and put $R_j(u)=e_ju$.
The normalized lifts
\[
u_j
=
R_{j+1}\circ p_{\omega_j}\circ R_j,
\qquad 0\leq j<r,
\]
are increasing degree-one homeomorphisms
of $\mathbb R$.
Choose continuous signed lifts
$l_{n,\omega}$ of the cocycle
$\theta_{n,\omega}$.
Set $m=r-L$, and define
\[
D_0
=
R_r\circ l_{n,\omega_{r-1}}
\circ\cdots\circ l_{n,\omega_0},
\]
\[
U_m=u_{m-1}\circ\cdots\circ u_0,
\qquad
D_{\mathrm{rem}}=D_0\circ U_m^{-1}.
\]
We shall retain the first $m$
steps and correct only the remaining
$L$ steps.

Let $A_-$ and $A_+$ be the endpoint
products of
$u_m,\ldots,u_{r-1}$.
By the width estimate,
\[
A_+(\omega,u)-A_-(\omega,u)>4
\]
for all $\omega\in V_r$ and
$u\in\mathbb R$.
Since $V_r$ is zero-dimensional, we may choose
a locally constant function
$\kappa:V_r\to\mathbb Z$ such that
\[
1<
D_{\mathrm{rem}}(\omega,0)+\kappa(\omega)
-A_-(\omega,0)
<3
\]
for every $\omega\in V_r$.
Since the difference of two increasing
degree-one homeomorphisms has
oscillation at most $1$, this gives
\[
A_-(\omega,u)
<
D_{\mathrm{rem}}(\omega,u)
+\kappa(\omega)
<
A_+(\omega,u)
\]
for every $\omega$ and $u$.

Proposition~\ref{prop:endpoint-interpolation}
therefore gives continuous increasing
degree-one lifts
$b_j$, $m\leq j<r$, satisfying
\[
\sup_{\omega,u}|b_j(\omega,u)-u_j(\omega,u)|
<\delta
\]
and
\[
b_{r-1}\circ\cdots\circ b_m
=
D_{\mathrm{rem}}+\kappa.
\]
For $j<m$, put $b_j=u_j$.
The complete normalized product is
then
\[
b_{r-1}\circ\cdots\circ b_0
=
D_0+\kappa.
\]
Returning to the original coordinates, define
\[
q_j=R_{j+1}\circ b_j\circ R_j.
\]
Thus
\[
q_j=p_{\omega_j}
\qquad (j<m),
\]
while
\[
\sup_{\omega,u}
|q_j(\omega,u)-p_{\omega_j}(u)|
<\delta
\qquad (m\leq j<r).
\]
Moreover,
\[
q_{r-1}\circ\cdots\circ q_0
=
l_{n,\omega_{r-1}}\circ\cdots
\circ l_{n,\omega_0}
+e_r\kappa.
\]
The integer disappears upon
passing to $\T$, so the corrected return maps agree exactly.

Write $\bar q_j$ for the circle map
induced by $q_j$, and let
\[
Q_{j,\omega}
=
\bar q_{j-1}\circ\cdots\circ\bar q_0,
\qquad
\Theta_{j,\omega}
=
\theta_{n,\omega_{j-1}}\circ\cdots
\circ\theta_{n,\omega_0},
\]
where both products are the identity
when $j=0$.
Define
\[
c_{n,\omega_j}
=
Q_{j,\omega}\circ\Theta_{j,\omega}^{-1},
\qquad 0\leq j<r.
\]
The formal value at $j=r$ is also
the identity, because the complete
return maps agree. Hence these
definitions are compatible at the
return bases.

The finite clopen tower partition
therefore determines a continuous
base-preserving homeomorphism
\[
C_n(\omega,z)
=
(\omega,c_{n,\omega}(z)).
\]
Each fiber map $c_{n,\omega}$
has degree one. By
Lemma~\ref{lem:continuous-lifts},
$C_n$ is isotopic to the identity,
and its pullback and inverse pullback
act identically on both $K$-groups. By construction,
\[
c_{n,\omega_{j+1}}
\circ\theta_{n,\omega_j}
\circ c_{n,\omega_j}^{-1}
=
\bar q_j,
\qquad 0\leq j<r,
\]
including the return step.
Since the corrected fiber maps
are within $\delta$ of the target
fiber maps, assertion (1) follows.

At every return base, $C_n$ is the
identity. Furthermore, for
$0\leq j\leq m$, the corrected
prefix $Q_{j,\omega}$ agrees with
the corresponding target prefix.
Consequently,
\[
C_nS_n^j(z)=T_n^j(z),
\qquad
z\in V_r\times\T,\quad
0\leq j\leq r-L.
\]
Moreover,
\[
(C_nS_nC_n^{-1})_\omega
=
\varphi_\omega,
\qquad \omega\in\Omega\setminus E,
\]
where
\[
E=
\bigcup_r\bigcup_{j=r-L}^{r-1}Q_{r,j}.
\]
This proves assertion (2).

For every
$\mu\in M(\alpha)\cup M(\alpha_n)$,
invariance and the equality of
the tower levels give
\[
\mu(Q_{r,j})=\mu(V_r).
\]
It follows that
\[
\mu(E)=L\sum_r\mu(V_r),
\qquad
1=\sum_r r\mu(V_r).
\]
Therefore
\[
\mu(E)
=
\sum_r\frac{L}{r}\,r\mu(V_r)
<\eta.
\]
This proves assertion (3), uniformly
over all the invariant measures
under consideration. Assertion (4)
has already been established.

Choose sequences $\delta_k,\eta_k\to0$.
For each $k\geq1$, let $N_k$ be sufficiently
large that the preceding construction applies
with $\delta=\delta_k$ and $\eta=\eta_k$
whenever $n\geq N_k$.
We may further assume that
\[
N_{k+1}>N_k\geq k
\]
and that $N_k$ is larger than any prescribed
lower bound.

For $N_k\leq n<N_{k+1}$, let $C_n$
be a homeomorphism obtained from the
construction with parameters
$(\delta_k,\eta_k)$.
Set $C_n=\mathrm{id}$ for $n<N_1$.
Then
\[
\sup_{\omega\in\Omega}
d_{\Homeo}\bigl(
c_{n,\alpha_n\omega}\circ
\theta_{n,\omega}\circ
c_{n,\omega}^{-1},
\varphi_\omega
\bigr)
<\delta_k
\]
whenever $N_k\leq n<N_{k+1}$.
Since $\alpha_n\to\alpha$ uniformly
and $\delta_k\to0$, it follows that
\[
C_nS_nC_n^{-1}\longrightarrow\taf
\]
uniformly. By
Lemma~\ref{lem:inverse-convergence},
the inverses converge uniformly as well.

Finally, if $S_n=h_nGh_n^{-1}$, set
$k_n=C_nh_n$. Then
\[
k_nGk_n^{-1}
=
C_nS_nC_n^{-1}.
\]
Since $C_n^*$ acts identically on
both $K$-groups,
$k_n$ and $h_n$ induce the same
maps on these groups. The same
holds for their inverses.
This completes the proof.
\end{proof}

\begin{rem}\label{rem:return-obstruction}
Proposition~\ref{prop:nearby-base-repair}
shows that the fiberwise conjugacies can be
modified to obtain uniform dynamical
approximation without changing their
induced maps on $K$-theory.
The correction is implemented by
homeomorphisms isotopic to the identity,
and hence is invisible at the level of
$K$-theory. Thus the dynamical correction can be
carried out while retaining the
prescribed $K$-theoretic information.
\end{rem}

\section{Rigidity and Approximate Unitary Equivalence}\label{sec:rigid}
\subsection{Tracial rank zero and classification}

We first recall the structural consequences of
classical rigidity and the classification results
needed for our purpose.

\begin{prop}\label{prop:rigid-TR0}
Let $A$ and $B$ be the crossed products
associated with the minimal classically rigid
systems in Theorem~\ref{thm:main}.
Then $A$ and $B$ are unital simple separable
nuclear $C^*$-algebras satisfying the UCT, with
\[
\tr(A)=\tr(B)=0.
\]
Moreover, every ordered unital graded
$K$-isomorphism
\[
k=(k_0,k_1):K_*(A)\longrightarrow K_*(B)
\]
is induced by a unital $*$-isomorphism
$\Phi:A\longrightarrow B$.
\end{prop}

\begin{proof}
By \cite[Corollary~3.11]{LinMatui2005},
classical rigidity implies that
$\rho_A(K_0(A))$ and $\rho_B(K_0(B))$
are uniformly dense in
$\Aff(T(A))$ and $\Aff(T(B))$,
respectively.

Since the underlying spaces are
one-dimensional, the tracial rank zero
criterion of
\cite[Theorem~4.6]{LinPhillips2010}
applies. Consequently,
\[
\tr(A)=\tr(B)=0.
\]
The remaining structural properties follow
from minimality and the standard properties
of crossed products by $\mathbb Z$. The classification theorem
\cite[Theorem~5.2]{LinClassificationTR0}
then yields a unital $*$-isomorphism
$\Phi:A\to B$ inducing the prescribed
ordered graded $K$-isomorphism.
\end{proof}

\begin{lem}\label{lem:rigid-trace-separation}
For each of the classically rigid systems
under consideration, the tracial state
space of its crossed product is determined
by the pairing with $K_0$.

In particular, if $\Phi:A\to B$ induces
$k_0$, then the affine homeomorphism
\[
L:T(B)\longrightarrow T(A),
\qquad L(\tau)=\tau\circ\Phi,
\]
is uniquely determined by $k_0$.
It satisfies
\begin{equation}\label{eq:rigid-pairing}
\rho_A(x)(L(\tau))
=
\rho_B(k_0(x))(\tau)
\end{equation}
for all $x\in K_0(A)$ and
$\tau\in T(B)$.
\end{lem}

\begin{proof}
By the standard trace correspondence for
crossed products of free minimal systems,
tracial states are identified with invariant
Borel probability measures on the underlying
spaces.

Classical rigidity identifies these measures
with the invariant probabilities on the Cantor
bases. Since clopen subsets determine Borel
probability measures on a Cantor space,
the traces are separated by the classes of
the corresponding clopen projections in $K_0$.

The pairing identity follows from
$\Phi_{*0}=k_0$. If two affine trace maps
satisfy this identity, their values agree
on all clopen projections and hence on
all traces by rigidity. Thus $L$ is unique.
\end{proof}

\subsection{Uniform trace convergence under rigidity}
\begin{lem}\label{lem:uniform-rigid-lift}
Let $S:Z\to Z$ have a compact metrizable factor $(\Omega,\delta)$,
with projection $\pi$, and assume that $\pi_*$ is injective on $M(S)$.
For each $n$, let $\mathcal E_n$ be a family of pairs
$(\mu,\lambda)\in\operatorname{Prob}(Z)\times M(S)$.
Suppose that 
\[
\int(a\circ S-a)\,d\mu\longrightarrow0\quad(a\in C(Z)),\qquad
\int b\circ\pi\,d(\mu-\lambda)\longrightarrow0\quad(b\in C(\Omega)).
\]
uniformly on $\mathcal{E}_n$. Then $\int a\,d(\mu-\lambda)\to0$ uniformly on $\mathcal E_n$
for every $a\in C(Z)$.
\end{lem}

\begin{proof}
Otherwise, for some $a\in C(Z)$ and
$\varepsilon>0$, there are
$(\mu_n,\lambda_n)\in\mathcal E_n$
such that
\[
\left|\int_Z a\,d(\mu_n-\lambda_n)\right|
\geq\varepsilon.
\]
Passing to a weak-* convergent subsequence,
the hypotheses imply that the limiting
measures $\mu,\lambda$ belong to $M(S)$
and satisfy $\pi_*\mu=\pi_*\lambda$.
By injectivity of $\pi_*|_{M(S)}$,
we have $\mu=\lambda$, contradicting
the above inequality.
\end{proof}

\begin{prop}
\label{prop:rigid-trace-lift}
Let $\Phi:A\to B$ be a unital
$*$-isomorphism inducing the ordered
graded $K$-isomorphism $k$, and put
\[
L:T(B)\to T(A),\qquad
L(\tau)=\tau\circ\Phi.
\]
Suppose that homeomorphisms
$\gamma_n:Y\to X$ satisfy
\[
\gamma_n\tbt\gamma_n^{-1}
\longrightarrow\taf
\]
uniformly. Let
\[
\chi_n:C(X)\to C(Y),\qquad
\chi_n(f)=f\circ\gamma_n.
\]
Assume that, for every finitely generated
subgroup $H\subset K_0(C(X))$,
\[
k_0\circ(j_{\taf})_{*0}|_H
=
(j_{\tbt}\circ\chi_n)_{*0}|_H
\]
for all sufficiently large $n$.

Then, for every $f\in C(X)$,
\begin{equation}\label{eq:rigid-trace-forward}
\sup_{\tau\in T(B)}
\left|
L(\tau)(j_{\taf}(f))
-
\tau(j_{\tbt}(f\circ\gamma_n))
\right|
\longrightarrow0.
\end{equation}
The reverse assertion holds with the
two systems interchanged.
\end{prop}

\begin{proof}
For each $\tau\in T(B)$, let $\nu_\tau$
and $\lambda_\tau$ be the invariant
probability measures corresponding to
$\tau$ and $L(\tau)$, respectively.
Set
\[
\mu_{n,\tau}=(\gamma_n)_*\nu_\tau.
\]
Thus $\lambda_\tau\in M(\taf)$, whereas
$\mu_{n,\tau}$ need not be $\taf$-invariant.

We first show that $\mu_{n,\tau}$ is
asymptotically $\taf$-invariant, uniformly
over $\tau\in T(B)$.
For every $a\in C(X)$, the
$\tbt$-invariance of $\nu_\tau$ gives
\[
\begin{aligned}
&\left|
\int_X(a\circ\taf-a)\,d\mu_{n,\tau}
\right|\\
&\quad=
\left|
\int_Y
(a\circ\taf\circ\gamma_n
-a\circ\gamma_n)\,d\nu_\tau
\right|\\
&\quad\le
\left\|
a\circ\taf\circ\gamma_n
-a\circ\gamma_n\circ\tbt
\right\|_\infty\\
&\quad=
\left\|
a\circ\taf
-a\circ\gamma_n\tbt\gamma_n^{-1}
\right\|_\infty
\longrightarrow0.
\end{aligned}
\]
Since the right-hand side is independent of $\tau$, the convergence is uniform over $T(B)$.

Now we let $U\subset\Omega_1$ be clopen, and
write
\[
p_U=1_{U\times\T}\in C(X).
\]
By the assumption, for all sufficiently
large $n$,
\[
k_0\bigl((j_{\taf})_{*0}[p_U]\bigr)
=
(j_{\tbt}\circ\chi_n)_{*0}[p_U].
\]
Using \eqref{eq:rigid-pairing}, we obtain
\[
\begin{aligned}
\lambda_\tau(U\times\T)
&=L(\tau)(j_{\taf}(p_U))\\
&=\rho_B\bigl(
k_0((j_{\taf})_{*0}[p_U])
\bigr)(\tau)\\
&=\tau(j_{\tbt}(p_U\circ\gamma_n))\\
&=\mu_{n,\tau}(U\times\T).
\end{aligned}
\]
This equality holds for every
$\tau\in T(B)$ and all sufficiently large $n$
depending only on $U$.

Now let $b\in C(\Omega_1)$ and
$\varepsilon>0$.
Choose a finite clopen partition
$\{U_1,\ldots,U_m\}$ of $\Omega_1$
and a step function
\[
s=\sum_{i=1}^m c_i1_{U_i}
\]
such that $\|b-s\|_\infty<\varepsilon$.
For all sufficiently large $n$, the
previous equality holds 
for each of $U_1,\ldots,U_m$. Hence
\[
\int_X(s\circ\pi)\,
d(\mu_{n,\tau}-\lambda_\tau)=0
\]
for every $\tau\in T(B)$.
Since both measures are probabilities,
\[
\sup_{\tau\in T(B)}
\left|
\int_X(b\circ\pi)\,
d(\mu_{n,\tau}-\lambda_\tau)
\right|
\le 2\varepsilon.
\]
Since $\varepsilon>0$ is arbitrary, we obtain
\[
\sup_{\tau\in T(B)}
\left|
\int_X(b\circ\pi)\,
d(\mu_{n,\tau}-\lambda_\tau)
\right|
\longrightarrow0
\qquad (b\in C(\Omega_1)).
\]
Now we may apply
Lemma~\ref{lem:uniform-rigid-lift}
to the families
\[
\mathcal E_n=
\{(\mu_{n,\tau},\lambda_\tau):
\tau\in T(B)\}.
\]
Classical rigidity of $\taf$
provides the required injectivity
of the projection on invariant
probability measures. Consequently,
for every $f\in C(X)$,
\[
\sup_{\tau\in T(B)}
\left|
\int_X f\,d(\mu_{n,\tau}-\lambda_\tau)
\right|\longrightarrow0.
\]
By the definitions of these measures,
this is precisely
\eqref{eq:rigid-trace-forward}.
The reverse assertion follows by
interchanging the two systems.
\end{proof}

\subsection{Approximate unitary equivalence of coefficient embeddings}
Write $\underline K(D)$ for ordinary $K$-theory together with all
$K_i(D;\Z/m)$, $m\ge2$, in both degrees.

\begin{lem}\label{lem:finite-total-K}
Let $\Omega$ be the Cantor space and $D=C(\Omega\times\T)$. Let $\phi,\phi_n:D\to E$ be
homomorphisms. If their induced maps agree eventually on every
finitely generated subgroup of $K_i(D)\,(i=0,1)$, then they agree eventually on every
finite subset of $\underline K(D)$.
\end{lem}

\begin{proof}
Since
\[
K_i(D)\cong C(\Omega,\Z),
\qquad i=0,1,
\]
both ordinary $K$-groups are torsion-free.
By the Bockstein exact sequence
\cite[Section~2.2]{MatuiHomomorphisms},
\[
K_i(D;\Z/m)\cong K_i(D)/mK_i(D).
\]
Thus every finite subset of
$\underline K(D)$ is determined by
finitely many ordinary $K$-classes.
The assertion follows from the
assumed eventual agreement.
\end{proof}

\begin{thm}[{\cite[Theorem~4.7]{MatuiHomomorphisms}; see also \cite[Definition~2.2]{MatuiHomomorphisms}}]
\label{thm:fixed-embedding-uniqueness}
Let $D$ be unital simple separable infinite dimensional with
$\tr(D)=0$, and fix a unital monomorphism 
\[a:C(W)\to D,\]
where $W$
is compact metrizable. Given finite $\mathcal F\subset C(W)$ and
$\epsilon>0$, there are finite $\mathcal P\subset\underline K(C(W))$,
finite $\mathcal H\subset C(W)$, and $\delta>0$ such that every
unital monomorphism $b:C(W)\to D$ satisfying
\[
a_*|_{\mathcal P}=b_*|_{\mathcal P},\qquad
\sup_{\tau\in T(D)}|\tau(a(h)) -\tau(b(h))|<\delta
\quad(h\in\mathcal H)
\]
admits a unitary $u\in D$ with
$\|ua(f)u^*-b(f)\|<\epsilon$ for $f\in\mathcal F$.
\end{thm}

This is the specialization of
\cite[Theorem~4.7]{MatuiHomomorphisms}
to unital monomorphisms.
Indeed, $D$ belongs to the class
$\mathcal T$ of
\cite[Definition~2.2]{MatuiHomomorphisms},
and the approximate multiplicativity
requirement is automatic.

\begin{prop}\label{prop:fixed-reference-realization}
Let $k:K_*(A)\to K_*(B)$ be an ordered
unital graded $K$-isomorphism, and fix
an isomorphism $\Phi:A\to B$ inducing $k$.
Suppose coefficient isomorphisms
$\chi_n,\eta_n$ satisfy the finite coefficient $K$-condition and
\eqref{eq:rigid-trace-forward}, together with its reverse for
$L^{-1}$, where $L(\tau)=\tau\Phi$.
Then there are unitaries $u_n\in B,w_n\in A$ with
\begin{align*}
\|\Ad(u_n)\Phi(j_{\taf}(f))-j_{\tbt}\chi_n(f)\|&\longrightarrow0,\\
\|\Ad(w_n)\Phi^{-1}(j_{\tbt}(g))-j_{\taf}\eta_n(g)\|&\longrightarrow0
\end{align*}
for every $f\in C(X)$ and $g\in C(Y)$.
\end{prop}

\begin{proof}
Set
\[
a=\Phi j_{\taf}:C(X)\to B,
\qquad
b_n=j_{\tbt}\chi_n:C(X)\to B.
\]
These are unital monomorphisms, with $a$
independent of $n$. By
Proposition~\ref{prop:rigid-TR0}, $B$ has tracial rank zero.

The finite coefficient $K$-condition and
Lemma~\ref{lem:finite-total-K} imply that
$a$ and $b_n$ agree eventually on every
finite subset of $\underline K(C(X))$.
Moreover, the assumed trace convergence gives
\[
\sup_{\tau\in T(B)}
|\tau(a(f))-\tau(b_n(f))|
\longrightarrow0
\qquad (f\in C(X)).
\]

Choose a dense sequence
$\{f_j\}_{j=1}^\infty$ in the unit ball
of $C(X)$, and put
$\mathcal F_m=\{f_1,\ldots,f_m\}$.
Apply Theorem~\ref{thm:fixed-embedding-uniqueness}
to the embedding $a$, with
$\mathcal F_m$ and $\epsilon=1/m$.
For each $m$, let
$\mathcal P_m\subset\underline K(C(X))$,
$\mathcal H_m\subset C(X)$, and
$\delta_m>0$ be given by
Theorem~\ref{thm:fixed-embedding-uniqueness}.
Choose strictly increasing integers
$N_m\geq m$ such that, for every $n\geq N_m$,
\[
a_*|_{\mathcal P_m}
=(b_n)_*|_{\mathcal P_m}
\]
and
\[
\sup_{\tau\in T(B)}
|\tau(a(h))-\tau(b_n(h))|
<\delta_m
\qquad(h\in\mathcal H_m).
\]

For $n\geq N_1$, set
\[
m(n)=\max\{m:N_m\leq n\}.
\]
Then $m(n)\to\infty$, and the uniqueness
theorem gives $u_n\in U(B)$ such that
\[
\|u_na(f)u_n^*-b_n(f)\|
<\frac1{m(n)}
\qquad(f\in\mathcal F_{m(n)}).
\]
For $n<N_1$, take $u_n=1$.
By density and contractivity, it follows that
\[
\|u_na(f)u_n^*-b_n(f)\|
\longrightarrow0
\qquad(f\in C(X)).
\]
Applying the same argument to the
embeddings
\[
\Phi^{-1}j_{\tbt},\quad j_{\taf}\eta_n:C(Y)\to A
\]
and using the reverse
$K$-group and trace conditions, give
unitaries $w_n\in U(A)$ and the second
asserted limit.
\end{proof}

\section{Proof of the main theorem}
\label{sec:proof}

\begin{proof}[Proof of Theorem~\ref{thm:main}]
Suppose first that condition~\textup{(1)}
holds, and choose sequences
$\sigma_n,\gamma_n,\Phi_n,\Psi_n$
as in Definition~\ref{def:approx-K-conjugacy}.
Set
\[
\chi_n(f)=f\circ\gamma_n,
\qquad
\eta_n(g)=g\circ\sigma_n.
\]
Since the isomorphisms $\Phi_n$ have
the same $KL$-class, they induce a
fixed ordered unital graded
$K$-isomorphism $k=(k_0,k_1)$.
Similarly, the maps induced by
$\Psi_n$ are $k_i^{-1}$.

For every $i=0,1$ and
$x\in K_i(C(X))$, the coefficient
norm convergence implies
\[
k_i((j_{\taf})_{*i}(x))
=
(j_{\tbt}\chi_n)_{*i}(x)
\]
for all sufficiently large $n$.
Indeed, sufficiently close projections
represent the same $K_0$-class,
and sufficiently close unitaries
represent the same $K_1$-class.
The same argument gives
\[
k_i^{-1}((j_{\tbt})_{*i}(y))
=
(j_{\taf}\eta_n)_{*i}(y)
\]
eventually for every
$y\in K_i(C(Y))$.
Choosing finite generating sets
and taking a common large $N$ proves
\eqref{eq:main-K-forward} and
\eqref{eq:main-K-backward}.
Thus condition~\textup{(2)} holds.

Conversely, assume condition~\textup{(2)}.
By Proposition~\ref{prop:simultaneous-normalization},
there are homeomorphisms
\[
H_n:Y\to X,\qquad K_n:X\to Y,
\]
with base maps $h_n$ and $\ell_n$, such that
\[
h_n\beta h_n^{-1}\longrightarrow\alpha,
\qquad
\ell_n\alpha\ell_n^{-1}\longrightarrow\beta
\]
uniformly. Moreover, for all sufficiently large $n$,
\[
\deg(H_n\tbt H_n^{-1})_\omega=s(\omega),
\qquad
\deg(K_n\taf K_n^{-1})_\eta=t(\eta)
\]
for every $\omega\in\Omega_1$ and
$\eta\in\Omega_2$, where $s$ and $t$
are the fiber orientation signs of
$\taf$ and $\tbt$, respectively.
Finally, for $i=0,1$,
\[
(j_{\tbt}H_n^*)_{*i}
=(j_{\tbt}\chi_n)_{*i},
\qquad
(j_{\taf}K_n^*)_{*i}
=(j_{\taf}\eta_n)_{*i}.
\]
Put
\[
S_n=H_n\tbt H_n^{-1}.
\]
Its base map is
$h_n\beta h_n^{-1}$, which converges
uniformly to $\alpha$, and its fiber
degrees eventually agree with those
of $\taf$, as shown before.
Since $\taf$ is minimal,
Proposition~\ref{prop:nearby-base-repair}
provides fiber-preserving
homeomorphisms $C_n:X\to X$,
isotopic to the identity, such that
\[
\gamma_n=C_nH_n
\]
satisfies
\[
d_{\Homeo(X)}
(\gamma_n\tbt\gamma_n^{-1},\taf)
\longrightarrow0.
\]

Applying the same proposition with
the systems interchanged gives
homeomorphisms $D_n:Y\to Y$,
isotopic to the identity, such that
\[
\sigma_n=D_nK_n
\]
satisfies
\[
d_{\Homeo(Y)}
(\sigma_n\taf\sigma_n^{-1},\tbt)
\longrightarrow0.
\]
Since $C_n$ and $D_n$ are isotopic
to the identity, their pullbacks
induce the identity maps on
$K_i(C(X))$ and $K_i(C(Y))$,
respectively, for $i=0,1$.
Consequently,
\[
(j_{\tbt}\gamma_n^*)_{*i}
=
(j_{\tbt}\chi_n)_{*i},
\qquad
(j_{\taf}\sigma_n^*)_{*i}
=
(j_{\taf}\eta_n)_{*i},
\quad i=0,1.
\]
In particular, the induced isomorphisms
\[
\widehat\chi_n(f)=f\circ\gamma_n,
\qquad
\widehat\eta_n(g)=g\circ\sigma_n
\]
still satisfy the finite coefficient
$K$-condition with respect to $k$.

By Proposition~\ref{prop:rigid-TR0},
there is a fixed isomorphism
$\Phi:A\to B$ inducing $k$.
Let
\[
L:T(B)\to T(A),
\qquad L(\tau)=\tau\circ\Phi.
\]
Since both systems are classically
rigid, Proposition~
\ref{prop:rigid-trace-lift},
applied in both directions, gives
\[
\sup_{\tau\in T(B)}
\left|
L(\tau)(j_{\taf}(f))
-
\tau(j_{\tbt}\widehat\chi_n(f))
\right|
\longrightarrow0
\]
for every $f\in C(X)$, and
\[
\sup_{\omega\in T(A)}
\left|
L^{-1}(\omega)(j_{\tbt}(g))
-
\omega(j_{\taf}\widehat\eta_n(g))
\right|
\longrightarrow0
\]
for every $g\in C(Y)$.
We may therefore apply
Proposition~
\ref{prop:fixed-reference-realization}
to obtain unitaries
$u_n\in U(B)$ and $w_n\in U(A)$
such that
\[
\begin{aligned}
&\|\Ad(u_n)\Phi(j_{\taf}(f))
-j_{\tbt}(f\circ\gamma_n)\|
\longrightarrow0,\\
&\|\Ad(w_n)\Phi^{-1}(j_{\tbt}(g))
-j_{\taf}(g\circ\sigma_n)\|
\longrightarrow0
\end{aligned}
\]
for every $f\in C(X)$ and
$g\in C(Y)$. Finally, we define
\[
\Phi_n=\Ad(u_n)\circ\Phi,
\qquad
\Psi_n=\Ad(w_n)\circ\Phi^{-1}.
\]
Inner automorphisms induce the
identity in $KL$, so
\[
[\Phi_n]=[\Phi],
\qquad
[\Psi_n]=[\Phi]^{-1}
\]
for every $n$.

As a conclusion, the homeomorphisms
$\gamma_n,\sigma_n$ satisfy both
dynamical limits and both
coefficient norm limits in
Definition~\ref{def:approx-K-conjugacy}.
Together with the fixed inverse
$KL$-classes and the paired
Elliott invariant induced by $\Phi$,
this proves condition~\textup{(1)}
and the final assertion of the theorem.
\end{proof}

% Acknowledgements and funding information may be added here.

% Bibliography: add your \bibitem entries in a thebibliography environment,
% or use \bibliographystyle{amsplain} and \bibliography{references}.

\begin{thebibliography}{99}

\bibitem{Blackadar1998}
B.~Blackadar,
\textit{$K$-Theory for Operator Algebras},
2nd ed.,
Mathematical Sciences Research Institute
Publications, vol.~5,
Cambridge University Press, Cambridge, 1998.



\bibitem{GiordanoPutnamSkau1999}
T.~Giordano, I.~F.~Putnam, and C.~F.~Skau,
Full groups of Cantor minimal systems,
\textit{Israel J. Math.}
\textbf{111} (1999), 285--320.




\bibitem{GiordanoPutnamSkau2004}
T. Giordano, I. F. Putnam, and C. F. Skau,
Affable equivalence relations and orbit structure
of Cantor dynamical systems,
\textit{Ergodic Theory Dynam. Systems}
\textbf{24} (2004), no.~2, 441--475.





\bibitem{HermanPutnamSkau1992}
R.~H.~Herman, I.~F.~Putnam, and
C.~F.~Skau,
Ordered Bratteli diagrams, dimension groups
and topological dynamics,
\textit{Internat. J. Math.}
\textbf{3} (1992), no.~6, 827--864.




\bibitem{LinClassificationTR0}
H. Lin,
Classification of simple $C^*$-algebras of tracial topological rank zero,
\textit{Duke Math. J.} \textbf{125} (2004), no.~1, 91--119.


\bibitem{LinFurstenberg2008}
H.~Lin,
Furstenberg transformations and approximate
conjugacy,
\textit{Canad. J. Math.}
\textbf{60} (2008), no.~1, 189--207.


\bibitem{LinMatui2005}
H. Lin and H. Matui,
Minimal dynamical systems on the product of the
Cantor set and the circle,
\textit{Comm. Math. Phys.} \textbf{257} (2005),
no.~2, 425--471.

\bibitem{LinMatui2005AC}
H. Lin and H. Matui,
Minimal dynamical systems and approximate conjugacy,
\textit{Math. Ann.} \textbf{332} (2005),
no.~4, 795--822.

\bibitem{LinMatui2006}
H. Lin and H. Matui,
Minimal dynamical systems on the product of the
Cantor set and the circle. II,
\textit{Selecta Math. (N.S.)} \textbf{12} (2006),
no.~2, 199--239.



\bibitem{LinPhillips2010}
H. Lin and N. C. Phillips,
Crossed products by minimal homeomorphisms,
\textit{J. Reine Angew. Math.} \textbf{641} (2010), 95--122.








\bibitem{Matui2005FullGroups}
H. Matui,
Approximate conjugacy and full groups of Cantor minimal systems,
\textit{Publ. Res. Inst. Math. Sci.}
\textbf{41} (2005), no.~3, 695--722.



\bibitem{Matui2006FullGroups}
H.~Matui,
Some remarks on topological full groups
of Cantor minimal systems,
\textit{Internat. J. Math.}
\textbf{17} (2006), 231--251.


\bibitem{MatuiHomomorphisms}
H. Matui,
Classification of homomorphisms into simple
$\mathcal{Z}$-stable $C^*$-algebras,
\textit{J. Funct. Anal.} \textbf{260} (2011),
no.~3, 797--831.






\bibitem{PimsnerVoiculescu1980}
M. Pimsner and D. Voiculescu,
Exact sequences for $K$-groups and
Ext-groups of certain cross-product
$C^*$-algebras,
\textit{J. Operator Theory}
\textbf{4} (1980), no.~1, 93--118.


\bibitem{Putnam1989}
I.~F.~Putnam,
The $C^*$-algebras associated with minimal
homeomorphisms of the Cantor set,
\textit{Pacific J. Math.}
\textbf{136} (1989), no.~2, 329--353.



\bibitem{Putnam2010Survey}
I.~F.~Putnam,
Orbit equivalence of Cantor minimal systems:
A survey and a new proof,
\textit{Expo. Math.}
\textbf{28} (2010), no.~2, 101--131.

\bibitem{WeggeOlsen1993}
N.~E.~Wegge-Olsen,
\textit{$K$-Theory and $C^*$-Algebras:
A Friendly Approach},
Oxford University Press, Oxford, 1993.


\end{thebibliography}

\end{document}